\documentclass[11pt]{amsart}
\usepackage{color}
\usepackage{pstool}
\usepackage{tikz-cd}

\usepackage[latin1]{inputenc} 
\usepackage[french]{babel}
\usepackage{psfrag}
\usepackage[T1]{fontenc}
 \usepackage{lmodern}
\usepackage{mathrsfs}
\usepackage{graphicx,epsfig,amscd,graphics,xypic, subcaption}
\usepackage{amsmath,amssymb}
\usepackage{enumitem}
\usepackage{rotating}

\usepackage[margin=4cm]{geometry}

\newcommand{\C}{\mathbb{C}}

\newcommand{\R}{\mathbb{R}}

\newcommand{\M}{\mathcal{MD}_{g,n}}
\newcommand{\T}{\mathcal{TD}_{g,n}}

\usepackage{lipsum}

\newtheorem{theorem}{Theorème}

\newtheorem*{theorem*}{Theorem}   
\newtheorem*{lemma*}{Lemme}   
  
\newtheorem*{conjecture}{Conjecture}  
\newtheorem{proposition}[theorem]{Proposition}
\newtheorem{defi}{Definition}           
\newtheorem{question}{Question} 
\newtheorem{lemma}[theorem]{Lemme}
\newtheorem*{remarque}{Remarque}

\newtheorem*{ex*}{Example \ref{exJ2} (continued)}

\title[Surfaces de dilatation]{Une invitation aux surfaces de dilatation}

\author[Selim Ghazouani]{Selim Ghazouani}
\address{Mathematics Institute, University of Warwick, Coventry CV4 7AL, U.K. }
\email{s.ghazouani@warwick.ac.uk}

\begin{document}
\maketitle

\begin{abstract}
Ce texte est une introduction aux \textit{surfaces de dilatation}. On  tente d'exposer les aspects géométriques et dynamiques du sujet: les espaces de modules, les feuilletages directionnels et la dynamique du flot de Teichmüller.  
\end{abstract}

\section{Introduction}

Une \textit{surface de dilatation} est une surface munie d'une structure géométrique modelée sur le plan complexe $\mathbb{C}$ via le groupe des similitudes (de la forme $z \mapsto az +b$) dont la partie linéaire $a$ est réelle (strictement) positive et pour lesquelles on autorise des points singuliers (sans quoi ces dernières ne sauraient exister!). Nous tâchons de donner dans ce texte une introduction (pratiquement) indépendante à ces objets, de motiver les principales questions qui nous semblent d'importance et de les replacer dans un contexte plus général, à la fois de géométrie et de systèmes dynamiques. L'exemple typique de surface de dilatation est l'espace obtenu en considérant un polygone du plan euclidien et en recollant des paires de côtés parallèles, comme dans la figure ci-dessous:

\begin{figure}[!h]
  \centering
  
  \includegraphics[scale=0.3]{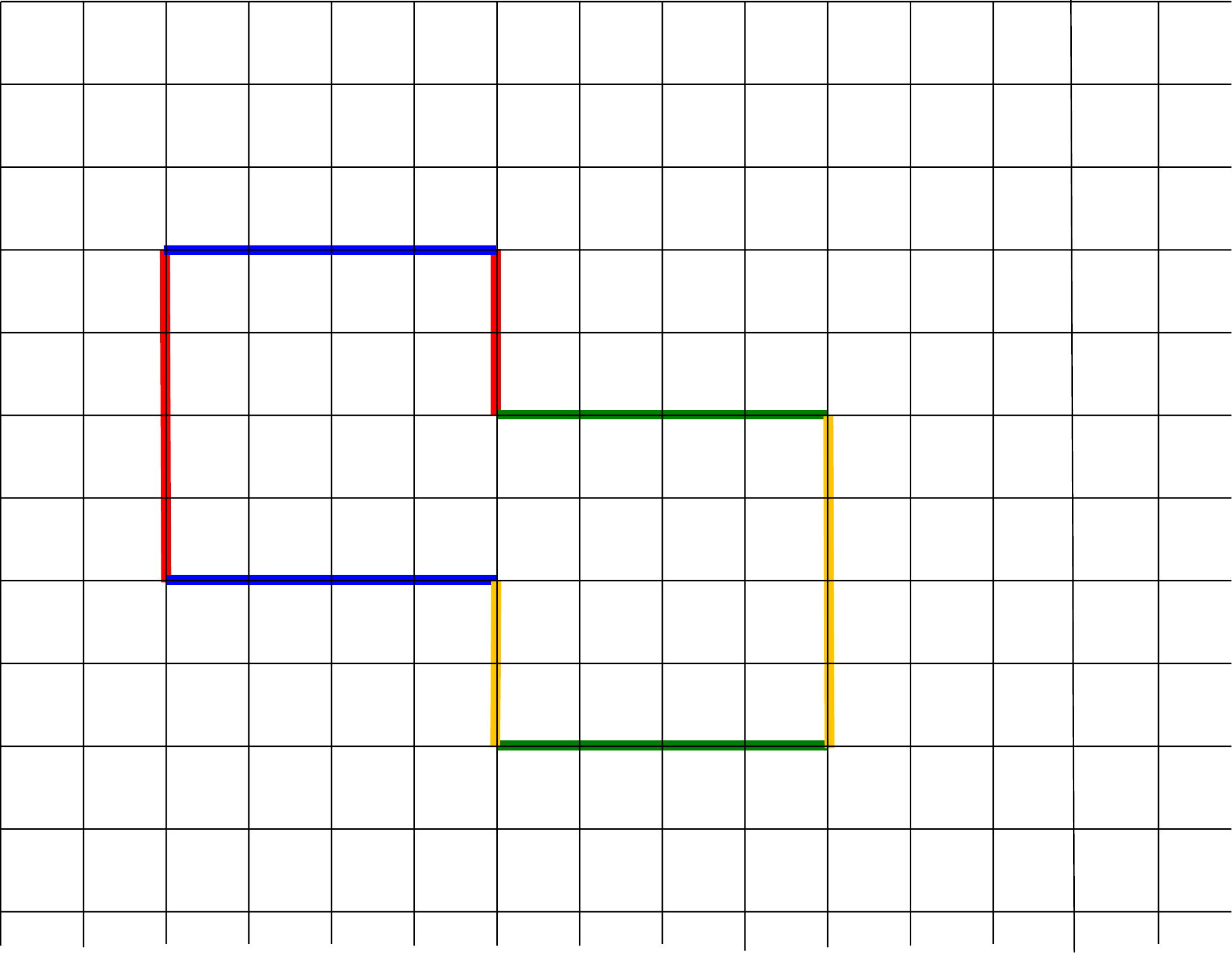}
  \caption{Un patron pour une surface de dilatation de genre $2$}
   \label{deuxchambres}
\end{figure}

\vspace{2mm} \noindent D'un point de vue \textit{géométrique}, ces structures se situent à mi-chemin entre les structures \textit{projectives complexes}(voir par exemple \cite{Dumas}) et les structures de \textit{translation}(voir \cite{Zorichsurvey}). Formellement elles peuvent être pensées soit comme un cas dégénéré des premières ou comme une généralisation des secondes. En pratique, nous pensons que le monde des surfaces de dilatation est un hybride équilibré entre les mondes projectifs et translations: \begin{itemize}

\item elles partagent avec les surfaces de translation l'objet \textit{ligne droite} qui donne une saveur "euclidienne" à leur étude;

\item elles n'ont pas de métrique naturelle et leur groupe d'holonomie "dilate" et "contracte" ce qui les rapproche de ce point de vue des structures projectives. 

\item à l'instar à la fois des structures projectives et de translation, elles ont une structure conforme compatible, ce qui ouvre la voie à la construction de leurs espaces de modules. 

\end{itemize}

\vspace{2mm} \noindent Nous pensons que l'intérêt particulier qu'on peut trouver chez les surfaces de dilatation vient de leur nature \textit{dynamique}. Elles portent en effet des familles de \textit{feuilletages directionnels} transversalement affines et leurs espaces de modules sont muni d'un \textit{flot de Teichmüller} qui agit comme \textit{opérateur de renormalisation} pour ces feuilletages. La construction cet opérateur est en tout point similaire au cas des surfaces de translation. Nous tenons cependant à souligner une différence fondamentale entre les mondes translation et dilatation.

\begin{itemize}

\item Les surfaces de translation sont à rapprocher des \textit{échanges d'intervalles linéaires} et la théorie de la renormalisation dans ce cas (déjà bien développée) s'intéresse aux propriétés \textit{ergodiques et combinatoires} de systèmes dynamiques unidimensionnels. Les échanges d'intervalles linéaires sont des isométries par morceaux et leurs propriétés géométriques sont aussi simples que possible.  

\item Il est par ailleurs bien connu que la théorie ergodique et combinatoire des \textit{échanges d'intervalles généralisés} se réduit complètement à celle des échanges linéaires. Les questions intéressantes pour ces échanges généralisés sont celles de nature géométrique(en d'autres termes qui relève de la "dynamique lisse"). Les feuilletages des surfaces de dilatations (et leurs avatars les échanges d'intervalles affines) sont la manifestation la plus simple de tels échanges d'intervalles(ou de feuilletages lisses sur des surfaces) dont les propriétés géométriques ne sont pas triviales. L'opérateur de renormalisation dans ce cas entretient des liens étroits avec ces propriétés géométriques et nous pensons que cela en fait un objet d'intérêt tout particulier.

\end{itemize}

\vspace{3mm} 

\paragraph*{\bf Organisation du texte} Les sections \ref{secGeometrie} à \ref{secEx} peuvent être prises comme un cours introductif, elles sont élémentaires et ne supposent aucun pré-requis (ni technique ni culturel). Dans la \textbf{section \ref{secGeometrie}}, nous définissons rigoureusement les surfaces de dilatation, les objets géométriques élémentaires qui leurs sont associés et prouvons des résultats basiques. Nous passons ensuite dans la \textbf{section \ref{secDyn}} aux feuilletages directionnels des surfaces de dilatation que nous définissons et pour lesquels nous passons en revue les propriétés élémentaires ainsi que leurs liens avec la géométrie. La \textbf{section \ref{secEsp}} est dédiée à la construction des espaces de modules de surfaces de dilatation et nous présentons en \textbf{section \ref{secEx}} une batterie d'exemples basiques de surfaces de dilatation. 

\vspace{3mm} \'A partir de la section \ref{teichmuller}, nous essayons d'introduire le lecteur aux propriétés de l'opérateur de renormalisation, la géométrie des espaces de modules et comment ces derniers interagissent avec les propriétés des feuilletages directionnels des surfaces de dilatation. Cette partie du texte nécessite probablement une plus grande familiarité avec les espaces de modules de structures géométriques et la renormalisation en systèmes dynamiques, bien que nous ayons tenu à éviter les références parachutées à des résultats connus des spécialistes. Dans la \textbf{section \ref{teichmuller}}, on introduit l'opérateur de renormalisation appelé \textit{flot de Teichmüller} et on expose une analogie entre ce flot et l'action du flot géodésique sur une variété hyperbolique de volume infini. Nous illustrons un cas où cette analogie s'avère suffisamment robuste pour prouver des résultats non-triviaux sur le comportement générique d'une classe de feuilletages transversalement affines sur le tore. Nous discutons dans la \textbf{section \ref{secAct}} de l'action naturelle du groupe $\mathrm{Sl}_2(\mathbb{R})$ et enfin nous terminons en \textbf{section \ref{secLien}} par une discussion plus poussée sur le cadre de dynamique unidimensionnelle dans lequel les questions que nous étudions se placent.

\section{Géométrie des surfaces de dilatation}
\label{secGeometrie}

\subsection{Définitions et préliminaires}

Il nous semble que le langage des $(G,X)$-structures peut s'avérer pratique pour introduire certaines notions importantes pour l'étude des surfaces de dilatation. Ce formalisme n'est cependant pas nécessaire et toute cette discussion peut être conduite dans un langage plus élémentaire, parti pris que nous avons choisi.

\vspace{2mm}. 

\noindent Une structure de dilatation sur une surface $\Sigma$ est la donnée d'un ensemble fini $S \subset \Sigma$ et un atlas $(U,\varphi_U)$ à valeurs dans $\mathbb{C}$ sur $\Sigma \setminus S$ tel que

\begin{enumerate}
\item les changements de cartes sont de la forme $z \mapsto az +b$ avec $a$ un réel strictement positif et $b \in \mathbb{C}$;

\item la structure ainsi définie sur $\Sigma \setminus S$ s'étend à chaque point de $S$ en une singularité de type conique-dilatation.

\end{enumerate}

\noindent Nous définissons maintenant formellement les singularités de type \textit{conique-dilatation}. Il convient de garder en tête qu'une telle singularité est simplement ce qu'on obtient, aux sommets dans l'espace quotient, en recollant des paires de côtés parallèles d'un polygone. Le modèle le plus simple est celui obtenu en prenant $\C^*$ que l'on fend le long d'une demi-droite issue de $0$ et en recollant les deux demi-droites ainsi obtenues le long de la dilatation $z \mapsto \lambda z$. On obtient ainsi en $0$ un singularité de type \textit{conique-dilation} d'angle $2\pi$ et de facteur de dilatation $\lambda$. En utilisant plusieurs feuillets on construit de la même manière la singularité de type \textit{conique-dilation} d'angle $2k\pi$ et de facteur de dilatation $\lambda$ pour tout entier $k$ strictement positif. 

\noindent Plus formellement,  c'est une structure de dilatation sur le disque épointé qui se construit de la manière suivante. Le plan épointé $\mathbb{C}^*$ est muni de la structure euclidienne standard qu'on peut tirer en arrière à son revêtement universel $\tilde{\C^*}$. Soit $u$ un générateur de l'action du groupe de revêtement.  Pour tout réel strictement poisitif $\lambda$, il existe une unique application $f_{\lambda}$ qui relève à  $\tilde{\C^*}$ l'application $z \mapsto \lambda z$ et qui fixe les relevés à $\tilde{\C^*}$ des droites qui passent pas $0$. Pour tout $k>0$ $f_{\lambda} \circ u^k$ agit proprement et discontinuement sur $\tilde{\C^*}$ et le quotient est un disque épointé. Une singularité de type \textit{conique-dilatation} sur $\Sigma$ d'angle $2k\pi$ et de facteur de dilatation $\lambda$ est un point de $S$ dont un voisinage épointé dans $\Sigma$ est affinement équivalent à un voisinage épointé de ce quotient.

\vspace{2mm}\noindent Si $s \in S$ est un point singulier de la structure, on $k_s$ l'entier tel que l'angle en $s$ est $2\pi k_s$ et $\lambda_s>0$ le paramètre de dilatation quand on parcourt une courbe autour de $s$ dans le sens positif.

\begin{proposition}[Formule de Gauss-Bonnet]

Pour toute surface de dilatation de genre $g \geq 1$, on a

\begin{enumerate}
\item $\sum_{s \in S}{2k_s\pi} = 4\pi(2g-2) $

\item $\sum_{s \in S}{\log\lambda_s} = 0$
\end{enumerate}

\end{proposition}

\begin{proof}

Pour la première partie de la proposition, il suffit de trianguler la surface avec une triangulation géodésique dont l'ensemble des sommets contient(et seulement contient!) les singularités, ce qui est facile à faire. La démonstration est à partir de là identique à celle du cas des surfaces de translation (écrire la caractéristique d'Euler en fonction du nombre de faces, sommets et arêtes et manipuler pour faire apparaitre le défaut d'angles aux singularités). 

\noindent Le second point se prouve en faisant apparaitre au concept d'holonomie linéaire qu'on introduit deux paragraphes plus loin dans le texte. Le nombre $\sum_{s \in S}{\log\lambda_s}$ est l'holonomie linéaire d'une courbe triviale en homologie et doit donc s'annuler.

\end{proof}

\noindent On déduit en particulier de cette formule de Gauss-Bonnet que l'unique surface compacte qui peut porter une structure de dilatation sans singularité est le tore.

\subsection{Objets vivants sur les surfaces de dilatation}

Un principe général dans l'étude des structures géométriques sur les variétés est que les objets bien définis sur le modèle et invariants par l'action du groupe de structure restent bien définis sur les variétés modelées sur ce modèle via l'action de ce groupe de structure. Dans notre cas, le modèle est le plan complexe $\mathbb{C}$ et le groupe de structure est le groupe affine complexe dont les parties linéaires sont réelles qu'on note soit $\mathbb{R}^*_+ \ltimes \C$ ou $\mathrm{Aff}_{\R^*_+}(\C)$. 

\noindent Les objets de cette géométrie sont les suivants.

\begin{enumerate}

\item Les \textit{droites} et les \textit{segments de droite}. En effet, l'action de $\mathbb{R}^*_+ \ltimes \C$ sur $\C$ préserve les lignes droites et une courbe sur une surface est dite \textit{droite} ou \textit{géodésique} si son image par une carte est toujours un bout de droite ou de segment de droite dans le plan. 

\item Les \textit{directions}. L'angle que forme une droite avec l'horizontale est bien défini car invariant par le groupe de structure. C'est une conséquence du fait que les parties linéaires de $\mathbb{R}^*_+ \ltimes \C$ sont réelles positives. 

\item Une simple extension du cas précédent nous dit que les \textit{angles} entre paires de droites/segments est bien défini.

\end{enumerate}

Les principales pertes par rapport au cas des surfaces de translation sont les notions de distance et d'aire. C'est une différence majeure qui a des conséquences importantes pour les systèmes dynamiques naturellement associés à nos structures, comme nous le verrons dans la suite. Mais cette perte de structure est aussi un gain en flexibilité.

\subsection{Holonomie linéaire}

Comme mentionné au début de cette section, le langage des $(G,X)$-structures serait tout à fait adapté à la discussion qui va suivre et le lecteur familier avec ce formalisme reconnaitra sans peine les notions de développante et d'holonomie cachées dans ce qui va suivre. 

\vspace{2mm} \noindent On veut définir pour toute courbe fermée une quantité qui encode combien la structure est "dilatée" le long de cette courbe. Considérons une courbe fermée $\gamma$ tracée sur une surface de dilatation $\Sigma$ qui évite les points singuliers. En partant d'un point $p$, on peut choisir une carte de la structure de dilatation en $p$ et chercher à l'étendre le long de $\gamma$ autant que faire se peut. Il n'y a pas d'obstruction à cette extension jusqu'à ce qu'on revienne près de $p$. On a défini par ce procédé deux cartes près de $p$ celle du "début" et celle de la "fin". Pour étendre la carte à tout $\gamma$ il faut que ces deux cartes coïncident. Ces deux cartes, par définition, diffèrent d'une transformation affine de la forme $z \mapsto \rho z + b$ avec $\rho$ un réel strictement positif. Le lecteur vérifiera sans peine que 

\begin{itemize}

\item ce nombre $\rho$ est indépendant du choix initial de la carte près de $p$;

\item indépendant du choix du point $p$;

\item indépendant du choix d'un représentant la classe d'homotopie libre de $\gamma$.

\end{itemize}

\noindent On définit ainsi pour toute courbe $\gamma$ un nombre $\rho(\gamma)$. On vérifie sans trop de difficulté que $\rho$ définit un morphisme de groupe

$$\rho : \mathrm{H}_1(\Sigma^*, \mathbb{Z}) \longrightarrow \mathbb{R}^* $$ où $\Sigma^*$ est la surface $\Sigma$ épointée des points singuliers de sa structure de dilatation. On appelle ce morphisme \textit{l'holonomie linéaire} de la surface $\Sigma$.

\subsection{Tores de Hopf et cylindres}

Les \textbf{tores de Hopf} sont probablement de bons candidats pour le titre de "surfaces de dilatation les plus élémentaires". Il s'agit du quotient du plan complexe privé de zéro par l'action d'une dilatation fixant zéro. Il s'agit d'une surface difféomorphe au tore $\mathbb{T}^2$ sans singularités. 

\vspace{2mm}

\noindent Une manière géométrique de représenter une telle surface est de considérer un anneau, la zone situé entre deux cercles concentriques de rayons respectifs $r$ et $R$. Le tore de Hopf est obtenu en recollant ces deux cercles via la dilatation de rapport $\lambda = \frac{R}{r}$. L'espace quotient est \textit{le tore de Hopf} de multiplicateur $\lambda$. 

\begin{figure}[!h]
  \centering
  \includegraphics[scale=0.4]{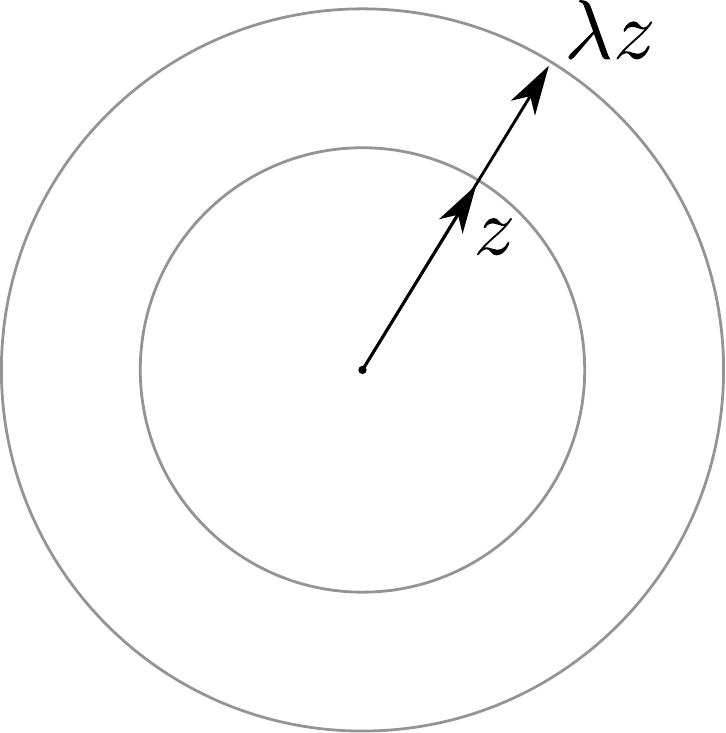}
  \caption{Un patron pour un tore de Hopf}
\end{figure}

\paragraph*{\bf Cylindres} Nous passons maintenant à la construction de ce que nous appelons \textit{cylindres de dilatation}. Nous les construisons de manière tout à fait analogue aux tores Hopf. La seule différence est qu'au lieu de considérer le plan complexe entier, nous nous restreignons à des secteurs angulaires. 
\noindent Le \textit{cylindre} $\mathcal{C}_{\theta, \rho}$ \textit{d'angle $\theta$ et de multiplicateur $\rho$} est le quotient d'un secteur angulaire d'angle $\theta$ centré en $0$ par l'action de la dilatation de facteur $\rho > 1$ qui fixe $0$. 

\vspace{1mm}
\noindent Dans la suite, quand nous disons qu'une surface de dilatation \textit{contient un cylindre}, nous entendons par là qu'il y a un plongement "de dilatation" (\textit{i.e.} qui préserve la structure de dilatation) du cylindre dans cette surface.

\subsection{Triangulations}

Dans la suite de ce texte, on entend par triangulation d'une surface de dilatation une triangulation dont l'ensemble des sommets est exactement l'ensemble des singularités et dont les arêtes sont des segments géodésiques. L'existence de telles triangulations est conceptuellement importante car elle permet de représenter une surface triangulable comme un simple recollement de polygones. 

\noindent Si toute surface de translation admet une triangulation, ce n'est pas le cas de toutes les surfaces de dilatation. En effet, il existe une obstruction facile à mettre en évidence:

\begin{proposition}

Soit $\Sigma$ une surface de dilatation qui contient un cylindre d'angle au moins $\pi$. Alors $\Sigma$ n'est pas triangulable. 

\end{proposition}

\noindent Un remaquable théorème du à Veech nous dit que cette obstruction est en fait la seule:

\begin{theorem}[Veech]
\label{thmtriangulations}
Une surface de dilatation admet une triangulation si et seulement elle ne contient pas de cylindre d'angle plus grand que $\pi$.
\end{theorem}

\noindent On pourra trouver la preuve de ce théorème dans l'appendice à l'article \cite{DFG}.

\section{Dynamique des feuilletages directionnels}
\label{secDyn}
\noindent Dans ce qui suit, $\Sigma$ est une surface de dilatation quelconque.

\subsection{Feuilletages directionnels}

Fixons une direction $\theta \in S^1$. L'ensemble des droites en direction $\theta$ sur $\Sigma$ forme ce qu'on appelle le \textit{feuilletage directionnel} en direction $\theta$. Cette définition est exactement la même que celle des feuilletages directionnels des surfaces de translation et ne fait que les généraliser à cette classe d'objets. 

\vspace{2mm} \noindent Ainsi toute surface de dilatation porte une famille de feuilletages indexé par l'ensemble des directions $S^1$. Nous faisons quelques remarques sur ces feuilletages:

\begin{enumerate}

\item ces feuilletages sont orientables, de telle sorte qu'il fait sens de parler de \textit{futur} et de \textit{passé};

\item les feuilletages orientés par $\theta$ et $-\theta$ sont ensemblistiquement les mêmes, mais le passé de l'un est le futur de l'autre et réciproquement;

\item ces feuilletages sont singuliers aux points singuliers de la structure de dilatation, à ces points la singularité est une selle dont le nombre de branches est exactement le nombre de fois qu'il faut mettre $\pi$ pour obtenir l'angle de la singularité.

\end{enumerate}

\noindent Encore une fois, tout ceci ne fait que généraliser directement le cas des surfaces de translation. Nous discutons maintenant les nouveaux  comportements qui apparaissent. 

\vspace{3mm} \noindent Considérons tout d'abord un exemple très concret, le tore de Hopf, disons de multiplicateur $2$. Chacun de ses feuilletages directionnel présente le même comportement dynamique qui est facile à décrire. Fixons une direction quelconque $\theta$. La droite orientée par $\theta$ passant par $0$ se projette sur le tore de Hopf sur deux géodésiques fermées, qui sont des feuilles fermées du feuilletage directionnel associé.  

\vspace{2mm}

\noindent Regardons de plus près ce qui se passe près de ces feuilles fermées. Prenons un petit segment $I \simeq ]-\epsilon, \epsilon[$ transverse à une de ces feuilles et intéressons-nous à l'application de premier retour du feuilletage sur $I$. En fonction du choix de l'orientation du feuilletage, cette application est la multiplication ou la division par $2$. Disons que nous avons choisi le sens qui nous donne la division par $2$. Cela nous dit que toute feuille qui passe dans un voisinage suffisamment proche de cette feuille fermée est d'une part piégée dans son voisinage mais surtout que la feuille vient s'accumuler sur cette feuille fermée. Formellement, l'$\omega$-limite de toute feuille contenue dans un certain voisinage de notre feuille fermée est cette feuille fermée. Nous avons donc des \textit{orbites périodiques} qui agissent comme des \textit{attracteurs}. Pour rendre cette discussion tout à fait symétrique, il faut dire qu'une feuille fermée qui est un attracteur dans le \textit{passé} est un \textit{répulseur}. Une feuille fermée de cette qualité (l'application de premier retour du feuilletage est contractante ou dilatante) est appelée en toute généralité \textit{hyperbolique}. Il ne faut pas travailler beaucoup plus pour obtenir le résultat suivant

\begin{proposition}

Pour tout feuilletage directionnel d'un tore de Hopf, il existe exactement deux feuilles fermées $c^+$ et $c^-$. La feuille $c^+$ est l'$\omega$-limite de toute feuille qui n'est pas $c^-$ et $c^-$ est l'$\alpha$-limite de toute feuille qui n'est pas $c^+$.
\end{proposition}

\noindent En d'autre terme, la feuille $c^-$ est la \textit{"source"} de ce feuilletage et $c^+$ en est son \textit{"puit"}. Ce comportement diffère singulièrement du cas des surfaces de translation, l'application de premier retour d'un tel feuilletage étant toujours une translation, les feuilles fermées viennent avec un cylindre plat et sont en particulier \textit{instables}. 

\subsection{Cylindres et stabilité}

On poursuit notre analyse en donnant des définitions générales correspondant au comportement dynamique observé pour les tores de Hopf et on relie ces comportement à l'existence de cylindres. 

\begin{defi}
Un feuilletage directionnel est dit \textit{Morse-Smale} ou \textit{asymptotiquement périodique} si l'$\omega$-limite et l'$\alpha$-limite de toute feuille non-singulière est une feuille fermée.  
\end{defi}

Ce comportement \textit{asymptotiquement périodique} est étroitement lié à l'existence de cylindres. En effet, soit $\mathcal{C}$ un cylindre dans une surface de dilatation $\Sigma$. Pour toute direction dans le secteur angulaire défini par ce cylindre, le feuilletage directionnel associé contient une feuille fermée contenue dans $\mathcal{C}$ et cette feuille est hyperbolique. Cette situation est en fait la cas général:

\begin{proposition}
Une feuille fermée régulière d'un feuilletage directionnel est hyperbolique si et seulement si elle est contenue dans un cylindre de dilatation. 
\end{proposition}

\begin{proof}

Cela résulte de la rigidité locale des structures affines. Deux géodésiques fermées simples tracées sur deux surfaces différentes et qui ont même holonomie linéaire ont des voisinages isomorphes. On applique ce fait général à une surface de dilatation quelconque contenant une feuille fermée d'holonomie linéaire $\rho$ dans une certaine direction et un cylindre de dilatation de multiplicateur $\rho$ contenant une feuille fermée dans la même direction.
\end{proof}

\subsection{Feuilletages transversalement affines} 

Nous faisons un petit aparté consacré à la notion de {\it feuilletage transversalement affine}, dont les feuilletages de dilatation sont des cas particuliers. Les feuilletages que nous considérons sont toujours singuliers, orientables et les singularités sont toujours de type selle. 

\begin{defi}
Un feuilletage transversalement affine $\mathcal{F}$ d'une surface compacte orientable $\Sigma$ est un feuilletage muni d'une structure affine transverse. Plus formellement c'est la donnée d'une application continue sur le revêtement universel de $\Sigma$

$$ \mathrm{D} : \tilde{\Sigma} \longrightarrow \mathbb{R}$$ constante sur les relevés des feuilles de $\mathcal{F}$ qui est équivariante via l'action d'une représentation $$ \rho : \pi_1 \Sigma \longrightarrow \mathrm{Aff}^+_1(\mathbb{R}) \simeq \mathbb{R}^*_+\ltimes \mathbb{R}$$ appelée holonomie.

\end{defi}

\noindent De manière plus informelle (mais pas moins rigoureuse), c'est la donnée pour tout petit segment transverse au feuilletage d'une identification canonique de ce segment avec un intervalle de $\mathbb{R}$ telle que le transport parallèle le long du feuilletage induise des applications affine entre segments transverses. Cette notion généralise la notion de feuilletage mesuré, de la même manière que les surfaces de dilatation généralisent les surfaces de translation. Un feuilletage mesuré est le cas où la représentation d'holonomie est triviale.

\vspace{2mm}

\noindent Tout les feuilletages directionnels de surfaces de dilatation sont transversalement affines. Il y a plusieurs manières de se convaincre de ce fait, la plus directe étant probablement de vérifier la définition formelle. Si $\theta$ est une direction et $\Sigma$ une surface de dilatation, on peut prendre la développante de sa structure affine et la projeter sur la direction orthogonale à $\theta$. On vérifie alors que cette fonction définit bien une structure transversalement affine sur le feuilletage directionnel associé (c'est plus ou moins tautologique).

\vspace{3mm}

\noindent Il nous semble pertinent de remarquer que cette classe de feuilletage est très naturelle à regarder: elle est suffisamment rigide pour avoir des paramétrisations simples (c'est à peu de choses près une variété de dimension finie) mais suffisamment flexible pour "voir" une grande partie de la complexité des feuilletages sur les surfaces (ce que ne permet pas la classe des feuilletages mesurés qui restreint trop la dynamique topologique des feuilletages). Nous voulions aussi citer le résultat suivant, dû à Liousse, qui prédit le comportement générique d'un tel feuilletage.

\begin{theorem}[Liousse, \cite{Liousse}]
\label{liousse}
Dans l'espace des feuilletages transversalement affines, le sous-ensemble des feuilletages de type Morse-Smale forme un ouvert dense.
\end{theorem}

\noindent Avant d'ouvrir une discussion sur la notion de feuilletage générique, nous voulions indiquer le problème ouvert suivant:

\begin{question}
Est-ce que tout feuilletage transversalement affine dont les singularités sont de type selle peut-être réalisé comme feuilletage directionnel d'une surface de dilatation?
\end{question}

\subsection{Généricité au sens mesuré}

Il existe (au moins deux notions) de généricité pour des familles de dimension finie de systèmes dynamiques:

\begin{itemize}

\item la généricité \textit{topologique}; on dit qu'un comportement est générique si il est observé sur un $G_{\delta}$ dense de la famille;

\item la généricité \textit{mesurée} ou \textit{probabiliste}; un comportement est dans ce cas générique si il est observé sur un ensemble de mesure pleine par rapport à la mesure de Lebesgue.

\end{itemize}

D'aucuns pourraient argumenter que la notion la plus pertinente d'un point de vue physique est la seconde. En effet, un ouvert dense (comme dans le théorème de Liousse \ref{liousse}) peut avoir une mesure très petite. Dans le cas qu'on regarde, pn propose la conjecture suivante:

\begin{conjecture}
\label{conjMS}
Un feuilletage transversalement affine est génériquement {\bf au sens mesuré} Morse-Smale/asymptotiquement périodique.
\end{conjecture}

\noindent Cette conjecture est loin d'être totalement évidente, comme l'illustre l'exemple suivant des \textbf{langues d'Arnold}. On considère la famille suivante de difféormorphismes du cercle 

$$ r_{\alpha,\epsilon} : x \mapsto x + \alpha + \epsilon \sin(2\pi x) \ \mathrm{mod} \ 1 $$ 

\noindent Arnold prouve (on dit que la preuve peut être trouvée dans \cite{Arnold}) les deux faits suivants:

\begin{enumerate}
\item l'ensemble des paramètres $(\alpha, \epsilon)$ tels que $r_{\alpha,\epsilon}$ est Morse-Smale est un ouvert dense;

\item son complémentaire est de mesure strictement positive.
\end{enumerate}

\noindent Il existe donc des familles intéressantes de systèmes dynamiques pour lesquelles généricités topologique et mesurée ne coïncident pas. La conjecture \ref{conjMS} impliquerait en particulier que si nous avions remplacé la famille $(r_{\alpha,\epsilon})$ par une famille d'homéomorphismes du cercle affine par morceaux, la conclusion du théorème d'Arnold ne tiendrait plus au sens où le complémentaire des paramètres Morse-Smale devrait être de mesure nulle. 
\noindent Une famille d'homéomorphismes du cercle affines par morceaux est considéré dans \cite{Herman1}, cette famille est étudiée en plus grand détail dans \cite{CGTU}, article dans lequel les auteurs et autrices conjecturent que l'ensemble des systèmes non Morse-Smale est de mesure non nulle. On donne une preuve de l'invalidité de cette conjecture en section \ref{teichmuller}. On notera aussi qu'un théorème de nature similaire existe pour les homéomorphismes du cercle lisses par morceaux avec deux discontinuités de la dérivée, voir \cite{Khmelev}.  

\vspace{3mm}

\noindent La question de la généricité mesurée va être le fil conducteur du reste de notre exposition. Nous allons maintenant introduire rigoureusement les espaces des modules des surfaces de dilatation auxquels nous allons penser comme des espaces de paramètres de systèmes dynamiques. Nous voulons insister sur l'interaction entre géométrie et dynamique qui est le coeur de cette approche.

\section{Espace des modules de surfaces de dilatation}
\label{secEsp}

On définit dans cette section les espaces des modules de surfaces de dilatation. Ces espaces fournissent la paramétrisation naturelle et par conséquent le cadre dans lequel poser les questions de généricité dynamique discuté dans la section précédente.

\vspace{2mm} \noindent  Dans la suite du texte

\begin{itemize}
\item $g$ et $n$ sont des entiers tels que $2 - 2g - n < 0$;

\item $\Sigma_{g,n}$ est la surface topologique de genre $g$ avec $n$  points marqués qu'on note $p_1, \cdots, p_n$;

\item  $\lambda = (\lambda_1, \cdots, \lambda_n)$ est un $n$-uplet de nombres strictement positifs tels que $\prod{\lambda_i} = 1$;

\item $k = (k_1, \cdots, k_n)$ des entiers positifs(pouvant être nuls) tels que $\sum{k_i} = 2g-2$ 

\end{itemize} 

\vspace{2mm} \noindent On définit

$$ \mathcal{TD}^*_{g,n} = \big\{ \text{structures de dilatations }\Sigma_{g,n} \ \text{avec singularités aux points marqués} \big\}/_{isotopies} $$

$$ \mathcal{MD}^*_{g,n} = \big\{ \text{structures de dilatations}\Sigma_{g,n} \ \text{avec singularités aux points marqués} \big\}/_{difféomorphismes} $$

\paragraph{\bf Stratification de $ \mathcal{TD}^*_{g,n}$ et $\mathcal{MD}^*_{g,n}$} Les espaces $ \mathcal{TD}^*_{g,n}$ et $\mathcal{MD}^*_{g,n}$ peuvent être partitionnés en fixant le type des singularités (nous entendons par là en fixant l'angle et le paramètre de dilatation à chaque point singulier). Si $\lambda$ et $k$ sont comme ci-dessus on note simplement $\mathcal{TD}^*(k,\lambda)$ et $\mathcal{MD}^*(k,\lambda)$ les strates associées.

\vspace{2mm} \paragraph{\bf Surfaces triangulables} Il est important de remarquer que les espaces $\mathcal{TD}^*_{g,n}$ et $\mathcal{MD}^*_{g,n}$ contiennent un ouvert remarquable formés des surfaces qui sont \textbf{triangulable} ou de manière équivalente (d'après le théorème \ref{thmtriangulations}) qui ne contiennent pas de cylindre d'angle au moins $\pi$. Pour des raisons de nature dynamique, qui devraient s'éclaircir au fur et à mesure de la discussion, il semble que ce lieu est le plus intéressant à étudier. Cela motive la définition suivante:

\begin{equation*}
\begin{aligned}
\mathcal{TD}_{g,n} =  &\big\{ \text{structures de dilatation }\Sigma_{g,n} \ \text{avec singularités aux points marqués }   \\ 
 & \text{admettant une triangulation géodésique} \big\}/_{isotopies}
\end{aligned}
\end{equation*}  

\begin{equation*}
\begin{aligned}
\mathcal{MD}_{g,n} =  &\big\{ \text{structures de dilatation }\Sigma_{g,n} \ \text{avec singularités aux points marqués }   \\ 
 & \text{admettant une triangulation géodésique} \big\}/_{difféomorphismes}
\end{aligned}
\end{equation*}  

\noindent Plus ou moins par définition, $\mathcal{MD}_{g,n}$ (resp. $\mathcal{MD}^*_{g,n}$) est le quotient de $\mathcal{TD}_{g,n}$ (resp. $\mathcal{TD}^*_{g,n}$) par l'action du groupe modulaire $\mathrm{Mod}(g,n)$.

\vspace{2mm} \noindent De manière analogue, on définit les strates  $\mathcal{TD}(k,\lambda)$ et $\mathcal{MD}(k,\lambda)$. Ce sera souvent avec ces espaces que nous travaillerons. D'après un théorème de Veech, tous ces espaces sont des orbifoldes.

\begin{theorem}[Veech, \cite{Veech}]
 Les espaces $\mathcal{TD}_{g,n}$, $\mathcal{TD}^*_{g,n}$ et chaque $\mathcal{TD}(k,\lambda)$  ($\mathcal{MD}_{g,n}$, $\mathcal{MD}^*_{g,n}$ et $\mathcal{MD}(k,\lambda)$) peuvent être munis d'une structure de variété(resp. d'orbifolde) analytique naturelle. La dimension de $\mathcal{TD}_{g,n}$, $\mathcal{TD}^*_{g,n}$, $\mathcal{MD}_{g,n}$ et $\mathcal{MD}^*_{g,n}$ est  $ 6(g-1) + 3n$  et celle de $\mathcal{TD}(k,\lambda)$ et $\mathcal{MD}(k,\lambda)$ est $ 6(g-1) + 2(n+1)$.
\end{theorem}

\subsection{Feuilletage par parties linéaires}

\label{partielineaire}

Nous entamons maintenant une discussion sur les propriétés élémentaires de ces espaces de modules. Le premier aspect sur lequel il nous semble important de se pencher est l'existence d'un feuilletage \textit{"par parties linéaires"} (aussi appelé ailleurs dans la littérature feuilletage \textit{isoholonomique}). Ce feuilletage est localement défini par la prescription le long de chaque courbe fermée de la quantité de dilatation que voit le transport parallèle le long de cette courbe. Du point de vue des polygones, il s'agit simplement des rapports entre les longueurs des côtés qu'on recolle. Formellement, $\T$ est un ensemble de structures de dilatation avec un marquage et par conséquent la fonction suivante (dite d'\textit{holonomie linéaire})

$$
\begin{array}{ccccc}
\mathrm{H} & : & \T & \longrightarrow & \mathrm{H}^1(\Sigma_{g,n}, \R^*) \\
 & & \Sigma & \longmapsto & \rho(\Sigma) 
\end{array}
$$ est bien définie. Il s'agit d'une submersion (voir Veech \cite{Veech}, p.625 Theorem 7.4) et par conséquent ces lignes de niveau forment un feuilletage trivial de $\T$. Par ailleurs, comme l'application $\mathrm{H}$ est équivariante par rapport à l'action linéaire naturelle du groupe modulaire sur $\mathrm{H}^1(\Sigma_{g,n}, \R^*)$, ce feuilletage passe au quotient $\T / \mathrm{Mod}(\Sigma_{g,n}) = \M$. Nous appelons ce feuilletage \textit{feuilletage par parties linéaire} ou parfois \textit{feuilletage isoholonomique}.

\vspace{2mm} \noindent Une remarque importante est que cette discussion peut-être conduite de manière tout à fait identique en remplaçant $\T$ et $\M$ par $\mathcal{TD}(k,\lambda)$ et $\mathcal{MD}(k,\lambda)$ respectivement et nous définissons donc par extension le feuilletage par partie linéaire restreint aux strates de manière  similaire.

\subsection{Une action de $\mathrm{Sl}_2(\mathbb{R})$}
\label{action}

Dans ce paragraphe nous définissons une action de $\mathrm{Sl}_2(\mathbb{R})$ sur $\T$ and $\M$. Une fois de plus, cette action est juste l'extension de l'action standard pour les surfaces de translation au cas qui nous intéresse. Et une fois de plus, la discussion qui va suivre vaut aussi pour $\mathcal{TD}(k,\lambda)$ et $\mathcal{MD}(k,\lambda)$.

\vspace{3mm} Nous expliquons ici la manière formelle de définir cette action mais on peut aisément faire les choses de manière plus intuitive au niveau des polygones. 

\noindent On considère un atlas de dilatation $(U, \varphi_U)_{U \in \mathcal{U}}$ sur  $\Sigma_{g,n}$ où $\mathcal{U}$ est une collection d'ouvert de $\Sigma_{g,n}$ tel que pour tout $U \in \mathcal{U}$, $\varphi_U : U \rightarrow \C$ est un homéomorphisme définissant la structure de dilatation. Considérons maintenant $A \in \mathrm{Sl}_2\R$. On vérifie sans peine que  $(U, A \circ \varphi_U)_{U \in \mathcal{U}}$ défini aussi une structure de dilatation sur $\Sigma_{g,n}$. On a par conséquent bien défini une action de  $\mathrm{Sl}_2 \R$  sur $\mathcal{TD}^*_{g,n}$ et $\mathcal{MD}^*_{g,n}$. 

\vspace{2mm} \noindent Remarquons maintenant les propriétés suivantes.

\begin{itemize}
\item L'image d'une triangulation géodésique par l'action d'un élément de  $\mathrm{Sl}_2\R$ est encore une triangulation géodésique. Notre action préserve donc $\T \subset  \mathcal{TD}^*_{g,n}$ et $\M \subset \mathcal{MD}^*_{g,n}$ et définit par restriction une action sur $\T$ et $\M$. Dans le reste du texte, nous ne nous intéresserons essentiellement qu'à cette action.  

\item Cette action sur $\T$ et $\M$ est localement libre (ce qui n'est pas le cas sur $ \mathcal{TD}^*_{g,n}$ et $\mathcal{MD}^*_{g,n}$, voir \cite{DFG}). 

\item Cette action préserve les feuilles du feuilletage par parties linéaires. 

\item Cette action préserve la stratification par les $ \mathcal{TD}(k, \lambda)$ et $\mathcal{MD}(k, \lambda)$.

\end{itemize}

\paragraph{\bf Flot de Teichmüller} Nous finissons cette section en insistant sur le fait suivant: l'action de $\mathrm{Sl}_2\mathbb{R}$ restreinte au sous-groupe des matrices diagonales

 $\big\{ \begin{pmatrix}
e^t & 0 \\ 0 &e^{-t} 
\end{pmatrix} \ | \ t \in \R \big\}$ définit un flot sur tous nos espaces de modules qui est généralement appelé \textit{flot de Teichmüller}. Nous distinguons ce sous-groupe en particulier car sa dynamique est liée aux comportements dynamiques des feuilletages portés sur les surfaces de dilatation. Plus précisément, si on considère une surface donnée comme un point de l'espace des modules, son orbite sous l'action du flot de Teichmüller nous donne beaucoup d'informations le comportement dynamique de son feuilletage vertical. Le flot de Teichmüller est ce qu'on appelle souvent en dynamique un \textit{opérateur de renormalisation}. Nous tentons d'éclaircir tout ça dans la section \ref{teichmuller}.

\subsection{Dégénérescence, compactification et fonction $\Theta$  }

On cherche maintenant à caractériser les possibles dégénérescences de surfaces de dilatation. Cela nous semble être un problème légèrement différent du cas des surfaces de translation. En effet, ces dernières sont des objets essentiellement riemanniens et leurs dégénérescences sont bien comprises en utilisant des fonctions comme la \textit{systole} ou le \textit{diamètre}, qui reposent sur le fait qu'on peut mesurer des longueurs. On fait dégénérer des surfaces de translation en écrasant des courbes, en regroupant des points singuliers, ou en étirant des cylindres. Pour les surfaces de dilatation (et c'est de manière plus générale une caractéristique de la géométrie affine), il est difficile de dire que quelque chose qui vit sur la surface est petit/grand par rapport à la surface. 

\vspace{2mm}
\noindent On identifie une quantité intéressante qui permet de caractériser certaines dégénérescences. Si $\Sigma$ est une surface de dilatation, on définit

$$ \Theta(\Sigma) := \sup_{\mathcal{C} \ \text{cylindre plongé dans} \ \Sigma}{ \{ \text{angle de } \mathcal{C} \} }.$$ 

\noindent Autrement dit, $ \Theta(\Sigma)$ est le plus grand angle d'un cylindre affine plongé dans $\Sigma$. On a les propriétés suivantes 

\begin{itemize}

\item si $\Sigma$ est une surface de translation, $ \Theta(\Sigma) = 0$;

\item $\Sigma$ est triangulable si et seulement si $ \Theta(\Sigma) < \pi$.  

\end{itemize}

\noindent Un problème intéressant est de déterminer quelles sont les surfaces pour lesquelles  $ \Theta(\Sigma) = 0$. 

\begin{question}

Est-il vrai que $ \Theta(\Sigma) = 0$ si et seulement si $ \Sigma$ est de translation? De manière équivalente, est-ce que toute surface de dilatation dont l'holonomie linéaire est non-triviale contient un cylindre de dilatation?

\end{question}

\noindent On sait que la réponse à cette question est positive en genre $1$. Si on comprenait au moins de manière grossière les dégénérescences de surfaces de dilatation, on pourrait imaginer une stratégie inductive pour attaquer cette question. Plus généralement, la question suivante nous semble être d'un intérêt tout particulier.

\begin{question}

Existe-t-il des compactifications intéressantes de $\mathcal{MD}_{g,n}$? 

\end{question}

\noindent On liste ci-dessous les différentes manières de dégénérer(assez triviales) que nous connaissons 

\begin{enumerate}

\item $\Theta$ tend vers $0$;

\item $\Theta$ tend vers $\pi$;

\item il existe un cylindre d'angle minoré par une constante strictement positive dont le multiplicateur tend vers l'infini ou vers $1$;

\item deux(ou plus) points singuliers s'agrègent.

\end{enumerate}

\noindent La question de la caractérisation des dégénérescences semble être une direction de recherche intéressante.

\section{Exemples}
\label{secEx}

On décrit dans cette section des exemples particulier de surfaces de dilatation, des constructions systématique et dans certains cas nous incluons la description dynamique de leurs feuilletages directionnels. 

\subsection{En genre $1$} Contrairement au cas des surfaces de translation, le cas du genre $1$ présente déjà une certaine richesse. En effet, on peut introduire un nombre arbitraire de points singuliers, qui seront tous d'angle $2\pi$ mais dont les facteurs de dilatation seront non-triviaux.

\vspace{2mm} 

\noindent Nous commençons par régler le cas sans singularités.  Nous avons vu la construction des tores de Hopf qui sont les quotients de $\mathbb{C}^*$ par une dilatation $z \mapsto \lambda z$ avec $\lambda$ strictement positif et différent de $1$. Notons $T_{\lambda}$ le tore de Hopf correspondant à la dilatation de facteur $\lambda$. Nous introduisons maintenant une petite variante de cette construction: coupons $T_{\lambda}$ le long d'une géodésique fermée simple, disons en direction horizontale (ça ne change pas grand chose car les tores de Hopf sont complètement symétrique vis-à-vis des directions). On obtient alors un cylindre affine d'angle $2\pi$. Il a deux composantes de bord totalement géodésiques. Ce que nous allons faire est recoller ces deux composantes de bord de manière différente. Pour n'importe quel choix d'une paire de points, chacun appartenant à une composante de bord différente, il existe une unique manière de recoller les composantes de bord de manière affine tout en identifiant ces deux points. On peut fixer le premier point pour supprimer la redondance d'information, et on obtient que l'ensemble des recollements possible est en bijection avec une des composante de bord. Parmi toute ces manière de recoller, une seule correspond au tore de Hopf $T_{\lambda}$. On appelle ces surfaces \textit{tores de Hopf modifiés}.

\begin{proposition}

\begin{enumerate}
\item Tout tore de dilatation sans singularité est isomorphe à un revêtement cyclique d'un tore de Hopf modifié.

\item Les feuilletages directionnels de tout tore de dilatation sans singularité sons Morse-Smale, et ont exactement $2k$ feuille fermés simple ($k$ attractives et $k$ répulsives) où $k$ est le degré du revêtement cyclique du point précédent.
\end{enumerate}

\end{proposition}

\noindent Nous voulons poursuivre cette discussion en montrant qu'au royaume des tores de dilatation, on peut tout construire en recollant des cylindres de dilatations de manière très simple. Commençons par l'exemple avec deux singularités donné par le recollement de l'hexagone suivant:

\begin{figure}[!h]

\label{toredessin}  
  
  \centering
  \includegraphics[scale=0.3]{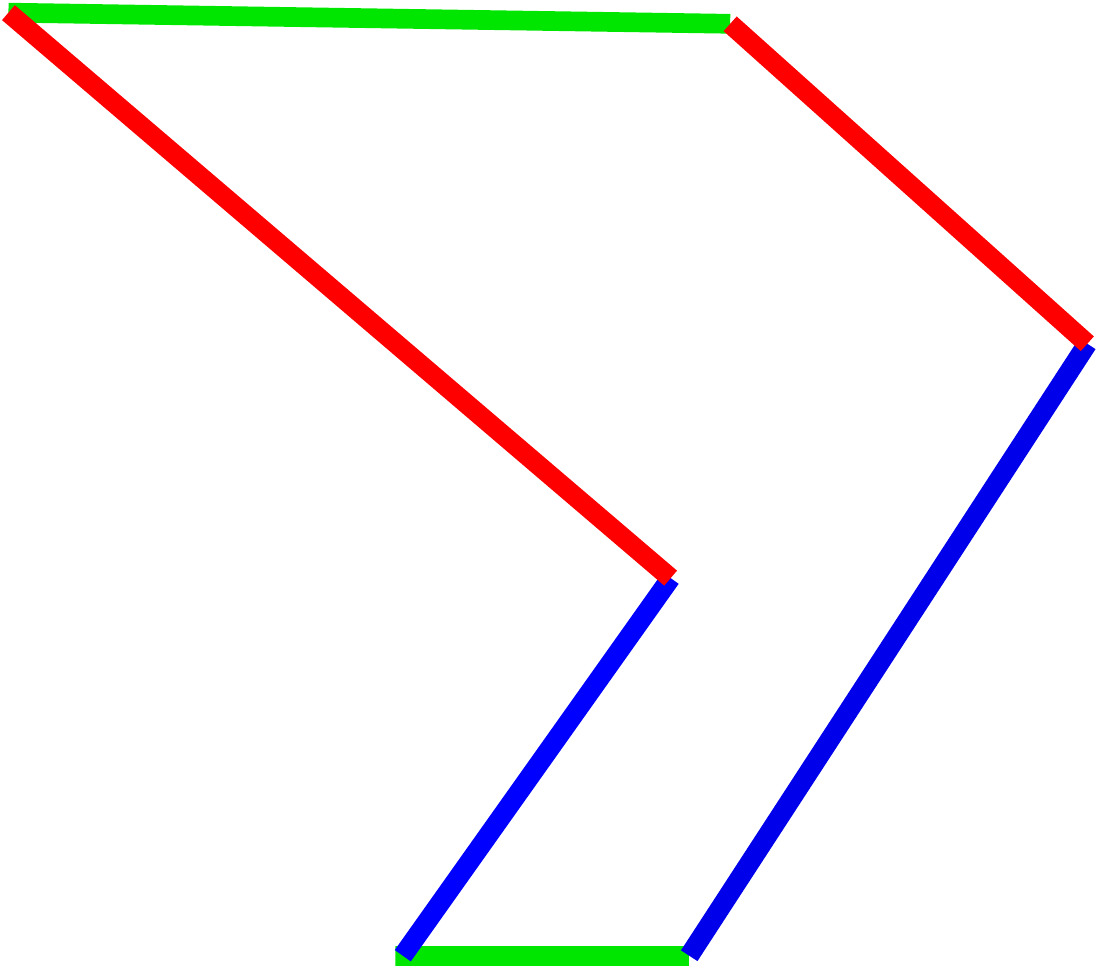}
  \caption{Un hexagone qu'on recolle pour former un tore de dilatation}
\end{figure}

\noindent On remarque que ce tore peut-être \textit{exactement} décomposé comme l'union de deux cylindres de dilatation de mêmes angles (qui sont sur la figure ci-dessus les enveloppes convexes respectives des paires de côtés rouges et bleus). Ce fait est assez général. On explique tout d'abord proprement la construction. On considère deux cylindres de dilatation (à bord totalement géodésique donc). On peut les recoller le long d'une composante de bord. Pour ce faire on est en règle général obligé de créer une singularité(ou plus d'ailleurs) d'angle $2\pi$ de facteur de dilatation qui est égal au rapport entre les facteurs de dilatation des cylindres qu'on recolle. On peut créer ainsi, en répétant l'opération, des cylindres de dilatation à bord géodésique avec en leur intérieur un nombre fini de singularités. Quand de surcroit, les deux composantes de bord sont dans la même direction, on peut les recoller pour refermer le cylindre et créer un tore.

\begin{remarque}
Le lecteur moins familier avec ce genre de constructions sera peut-être décontenancé par le manque de formalisme. Notre point de vue est que bien qu'il soit possible de rendre cette construction rigoureuse d'un point de vue formel, en introduisant les bonnes identifications entre les composantes de bord et en étendant la structure affine à l'espace quotient, il est préférable de ne pas alourdir l'exposition de formalisme décrivant des opérations qui sont au demeurant très intuitives. 
\end{remarque}

\noindent La proposition suivante établit le fait que tout tore de dilatation peut être construit de la façon décrite ci-dessus.

\begin{proposition}
\label{decomposition}

Soit $T$ un tore de dilatation avec $n$ singularités. Alors $T$ peut être décomposé comme l'union d'au plus $n$ cylindres de dilatation. 

\end{proposition}

\begin{proof}

La preuve repose essentiellement sur l'existence d'un cylindre(de dilatation ou plat) plongé dans $T$. Expliquons tout d'abord pourquoi cela implique l'existence de la décomposition. Un cylindre plongé peut être étendu jusqu'à ce que son bord rencontre une singularité (ou plusieurs singularités). Mais comme sur un tore de dilatation les singularités sont toutes d'angle $2\pi$, ce qu'il y a de l'autre côté du bord de ce cylindre est encore un cylindre car c'est un ouvert bordé par une géodésique fermé simple qui forme un angle $\pi$ à chaque singularité.  On feuillette ainsi le tore en géodésiques fermées simples qui doit se refermer une fois qu'on a épuisé toutes les singularités. 

\vspace{3mm}

On prouve maintenant l'existence d'un tel cylindre ou, de manière équivalente, l'existence d'une géodésique fermée (nécessairement simple). Choisissons une direction $\theta$ sur la surface, elle définit d'un point de vue topologique un feuilletage régulier orientable sur le tore. Soit ce feuilletage admet une feuille fermée et c'est gagné. Supposons donc qu'il n'en a pas. On peut dans ce cas aisément construire une courbe fermée simple transverse au feuilletage. La technique, standard en théorie des feuilletages, consiste à prendre un long rectangle bordé par deux feuilles du feuilletage dont une récurrente et d'utiliser une diagonale qu'on referme habilement en utilisant la récurrence. 

\vspace{2mm} \noindent Une fois muni d'une telle courbe fermée simple $\gamma$ transverse au feuilletage, on définit l'application de premier retour sur cette courbe. Nous affirmons qu'elle est bien définie \textit{i.e.} que toute feuille issue de cette courbe y revient. C'est une conséquence du fait que le complémentaire de cette courbe est un cylindre portant un feuilletage régulier, transverse au bord, sans feuille fermée. Un tel feuilletage est toujours trivial. 

\vspace{2mm} Ainsi la dynamique de notre feuilletage est-elle contenue dans celle de l'application de premier retour sur cette courbe fermée simple qui est donc un homéomorphisme du cercle $f_{\theta_0} : S^1 \longrightarrow S^1$. Celui-ci n'a pas de points fixes et est donc semi-conjugué à une rotation. Remarquons maintenant que pour des directions $\theta$ proches de $\theta_0$, $\gamma$ reste transverse au feuilletage directionnel et l'application de premier retour sur $\gamma$ qu'on note $f_{\theta}: S^1 \longrightarrow S^1$ et qui varie continument avec $\theta$. On peut dire plus (et c'est là qu'on utilise la géométrie): selon qu'on choisit $\theta$ plus grand ou plus petit que $\theta_0$ on a que $f_{\theta} = f_{\theta_0} + \epsilon(\theta)$ où $\epsilon$ est une petite perturbation qui est toujours \textbf{strictement positive} ou \textbf{strictement négative} et qui est \textbf{strictement croissante} avec $\theta$. Cela implique que le nombre de rotation en $\theta_0$ est strictement croissant. Parce que le nombre de rotation est une fonction continue, il va prendre des valeurs rationnelles. Or un homéomorphisme du cercle a un nombre de rotation rationnel si et seulement il a un orbite périodique. Il existe donc $\theta$ proche de $\theta_0$ tel que $f_{\theta}$ a une orbite périodique et par conséquent le feuilletage associé à une feuille fermée. Ceci achève la démonstration.
\end{proof}

\paragraph*{\bf Dynamique sur les tores de dilatations}

Le petit argument de la preuve de la proposition \ref{decomposition} nous montre que la question des possibles comportements dynamiques pour un feuilletage directionnel sur un tore peut-être ramené à l'étude d'homéomorphismes du cercle. 

\begin{itemize}
\item Si le feuilletage a une feuille fermée, on peut sans trop de peine que le feuilletage est \textbf{Morse-Smale}, à des cylindres plats près.

\item Si le feuilletage n'a pas de feuille fermée, il suffit de considérer une courbe transverse au feuilletage comme on l'a fait plus haut. L'étude dynamique du feuilletage se ramène à un homéomorphisme du cercle de nombre de rotation irrationnel. En utilisant la structure transversalement affine du feuilletage, on peut conjuguer cet homéomorphisme à un homéomorphisme affine par morceaux. Or on sait (voir \cite{Herman1}, Chapitre VI) qu'un tel homéomorphisme est \textbf{minimal} car conjugué à la rotation de même nombre de rotation. Cela vient du fait que la classe des  homéomorphismes affine par morceaux appartient à la classe P, ce qui revient à dire que leur dérivée est toujours à variation bornée.

\end{itemize}

\noindent En résumé, un feuilletage sur le tore est soit \textbf{minimal} (toutes les orbites sont dense) ou \textbf{asymptotiquement périodique}.

\subsection{En genre $2$}

Le cas du genre $2$ (et à plus forte raison le genre plus grand) présente une difficulté supplémentaire qui est l'existence de points singuliers dont l'angle est différent de $2\pi$. Ils sont responsables de la singularité des feuilletages associés et rendent les arguments combinatoires plus compliqués.

\vspace{2mm} \noindent Nous commençons par présenter une famille de surface que nous appelons \textit{"surfaces à deux chambres"}. Le comportement dynamique des feuilletages directionnels de ces surfaces a été récemment étudié dans \cite{Micheletmichel}.
\noindent Une \textit{"chambre"} est un tore de dilatation avec une composante de bord géodésique avec un point singulier sur le bord d'angle intérieur égal à $3\pi$. C'est exactement les surfaces qu'on obtient en recollant deux paires de cotés parallèles d'un pentagone:

\begin{figure}[!h]
  \centering
  \label{1chambre}
  \includegraphics[scale=0.4]{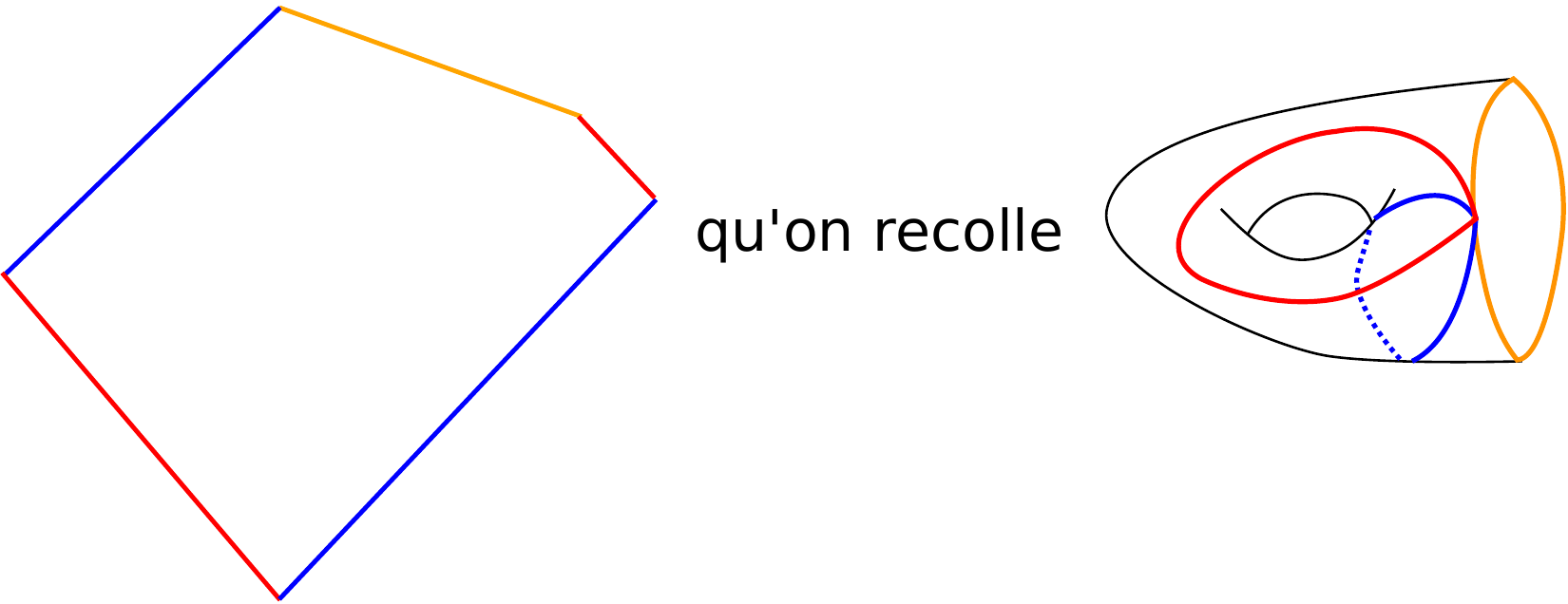}
  \caption{Une "chambre"}
\end{figure}

\noindent Une \textit{surface à deux chambres} est la surface qu'on obtient en recollant deux telles chambres le long de leurs composantes de bord en identifiant leur points singuliers. Les surfaces à deux chambres sont des surfaces de dilatations de genre $2$ avec une singularité d'angle $6\pi$ et de facteur de dilatation trivial. 

\noindent Ces surfaces ont de ça d'intéressant que la dynamique se scinde en deux parties: ce qui se passe dans chacune des chambres. Plus précisément, pour toute direction qui n'est pas celle de la "porte" de la chambre, toute feuille traverse la "porte" une unique fois et reste ultimement coincée dans une chambre dans le futur et dans l'autre dans le passé.

\noindent Nous référons aux articles \cite{DFG} et \cite{BowmanSanderson} pour plus de détails sur ce qui va être esquissé dans ce paragraphe. 

\vspace{3mm}

\noindent  Pour presque toute direction $\theta$ sur une de ces surfaces à deux chambres, le feuilletage possède exactement deux feuilles fermées hyperboliques, une dans chacune des chambres. La première(resp. la deuxième) est l'$\omega$-limite (resp. $\alpha$-limite) de toute feuille sauf de l'autre feuille hyperbolique.
\noindent Pour une surface donnée, l'ensemble des $\theta$ pour lesquels cet énoncé est vrai est un ouvert dense de mesure pleine.
\vspace{2mm}
\noindent  Par ailleurs, pour certaines directions qui forment un ensemble de Cantor de dimension de Hausdorff nulle, l'ensemble $\omega$-limite et/ou $\alpha$-limite des feuilles peut(peuvent) être un (des) ensembles transversalement Cantor (en d'autre terme une lamination).  Restreinte à ces laminations, le feuilletage possède la structure orbitale d'un feuilletage irrationnel sur le tore. Ces directions exceptionnelles sont essentiellement des \textit{flots de Cherry} (voir \cite{Cherry}).

\vspace{3mm}

Nous poursuivons par une remarque sur l'espace des modules de telles surfaces à deux chambres.

\begin{proposition}
L'ensemble des surfaces à deux chambres est connexe et forme une composante connexe de la strate à laquelle elles appartiennent.
\end{proposition}

La strate en question est $\mathcal{MD}\big( (2), (1) \big)$. C'est la strate de plus petite dimension en genre au moins $2$ et nous pensons qu'elle mérite une attention particulière, en tant que modèle pour le cas général. Nous posons la question suivante

\begin{question}
Combien la strate $\mathcal{MD}\big( (2), (1) \big)$ a-t-elle de composantes connexes?
\end{question}

\noindent Remarquons que la réponse à cette question est au moins $2$. En effet, il existe des surfaces qui ne sont certainement pas à deux chambres. C'est le cas par exemple des surfaces de translation, qui ont des déformations continues qui sont strictement de dilatation. Il suffit ainsi de prendre une surface de translation de genre $2$ avec une singularité et de considérer un voisinage dans $\mathcal{MD}\big((2),(1) \big)$ pour y trouver des surfaces qui ne peuvent être dans la composante des surfaces à deux chambres.

\section{Flot de Teichmüller et dynamique des feuilletages directionnels}
\label{teichmuller}

Nous rappelons le lecteur qu'une des principales motivations de l'étude des surfaces de dilatation est la manière dont elles interagissent avec la dynamique unidimensionnelle. On veut comprendre les feuilletages transversalement affine qu'elles portent et on espère que l'étude des espaces des modules va nous permettre de progresser dans cette direction. La question qui nous guide est la suivante 

\begin{center}

{ \it Quel est le comportement dynamique générique d'un feuilletage transversalement affine?}

\end{center}

\subsection{Renormalisation et philosophie} 

On explique dans ce paragraphe la philosophie sous-jacente au programme de recherche/méthodes que nous suggérons dans ce texte.   

\vspace{3mm}

\noindent Notre but est de démontrer que l'ensemble des feuilletages Morse-Smale est un ouvert dense de mesure pleine. On pense ici à un point dans un espace de module $\mathcal{MD}_{g,n}$ comme à un feuilletage en considérant le feuilletage vertical de la surface associée. Dans la suite, on note $\mathcal{P} \subset \mathcal{MD}_{g,n}$ l'ensemble des surfaces dont le feuilletage vertical est Morse-Smale/asymptotiquement périodique. 
\noindent Le fait que $\mathcal{P}$ est ouvert est un fait facile à vérifier qui découle de la stabilité des points fixes d'application contractante. La difficulté dans ce problème revient à estimer la "taille" de $\mathcal{P}$, tant d'un point de vue topologique que mesuré. 

\vspace{2mm}

\noindent On peut par exemple apercevoir "à l'oeil nu" des bouts de $\mathcal{P}$: sur une surface donnée, on voit des "gros cylindres" qui survivent à des déformations et permettent de localiser $\mathcal{P}$. Le problème vient du fait qu'il existe des orbites périodiques "très longues" et "topologiquement compliquées" qu'il est très difficile de distinguer de feuilletages minimaux (ou avec une dynamique plus complexe qu'asymptotiquement périodique). Par exemple, il peut arriver (on en donnera des exemples plus tard dans le texte) qu'un certain point dans $\mathcal{MD}_{g,n}$ corresponde à un feuilletage minimal. Dans un voisinage de ce point il y a plein de feuilletages asymptotiquement périodique qui approximent ce feuilletage de manière de plus en plus précise de telle sorte qu'à partir d'un moment il devient difficile de les différencier, à moins d'être capable de zoomer avec une précision infinie! 

\vspace{2mm} 
\noindent Nous avons introduit en \ref{action} de manière un peu sèche ce que nous avons appelé le \textit{flot de Teichmüller}. Il s'agit d'un flot sur $\mathcal{MD}_{g,n}$ qui réalise précisément le changement d'échelle dont nous avons besoin pour appréhender les feuilletages Morse-Smale qui sont tout de même très compliqués. Nous rappelons que le flot de Teichmüller $(g_t)$ réalise l'action des matrices diagonales

 $$ \begin{pmatrix}
e^t & 0 \\
0 & e^{-t} 
 \end{pmatrix}.$$ Ces dernières, quand $t$ est très grand, "rapetissent" le feuilletage vertical en comparaison aux autres directions, et permet de le remettre à l'échelle qui nous permet de le voir à l'oeil nu, et d'estimer la taille relative de $\mathcal{P}$ près de ce feuilletage. C'est pourquoi nous allons être particulièrement intéressé par la description de la dynamique de flot. Un tel opérateur en théorie des systèmes dynamiques est appelé \textit{opérateur de renormalisation}. 
 
\vspace{2mm} 
 
 Nous établirons dans les sections suivantes un dictionnaire(incomplet) entre les propriétés des orbites du flot de Teichmüller et les propriétés dynamiques des feuilletages transversalement affines portés par les surfaces de dilatation. Nous essayerons ensuite de comprendre comment la géométrie de $\mathcal{MD}_{g,n}$ peut contraindre la dynamique de $g_t$.

\subsection{Le flot de Teichmüller}

On introduit dans cette section des sous-ensembles de $\mathcal{MD}_{g,n}$ qui sont pertinents pour la compréhension du lien entre l'action du flot de Teichmüller et le problème du comportement générique des feuilletages des surfaces de dilatation. 

\begin{itemize}

\item $\mathcal{P}$ est l'ensemble des surfaces dont le feuilletage vertical est Morse-Smale.

\item $\mathcal{P}^*$ est l'ensemble des surfaces dont le feuilletage vertical contient une feuille fermée hyperbolique. 

\item $\mathcal{C}$ est le complémentaire dans $\mathcal{MD}_{g,n}$ de $\mathcal{P}^*$.

\item Si $\Sigma$ est une surface de dilatation, $\mathcal{C}_{\Sigma}$ est le sous-ensemble de $S^1$ tel que $r_{\theta} \cdot \Sigma \in \mathcal{C}$. 

\end{itemize}

\noindent Le résultat qu'on aimerait démontrer en toute généralité est le fait que $\mathcal{P}$ est de mesure pleine. Les feuilletages sur les surfaces de genre grand peuvent exhiber une complexité combinatoire qui rend l'analyse un peu pénible: on entend par là qu'il se passe des choses différentes dans différentes parties de la surface. Ce fait est responsable du fait que $\mathcal{P}$ et $\mathcal{P}^*$ sont des ensembles différents. Nous voulons dire par là que ce n'est qu'une difficulté technique et que nous concentrons nos efforts sur $\mathcal{P}^*$, le passage à $\mathcal{P}$ étant essentiellement une question d'architecture de preuve. 

\vspace{3mm}

\noindent On présente deux lemmes qui nous permettent de mieux comprendre le lien entre flot de Teichmüller et feuilletage sur les surfaces. L'idée générale est que l'ensemble $\mathcal{P}^*$ est un ouvert dense et son complémentaire $\mathcal{C}$ est un fermé qui ressemble à un ensemble de Cantor. Plus précisément, on regarde (ou plutôt on de réduire le problème à) la trace de $\mathcal{C}$ sur l'ensemble des directions de chaque surface, qu'on a décidé de noter $\mathcal{C}_{\Sigma}$. Il faut s'attendre à ce que ça soit un ensemble de Cantor et c'est la taille de cet ensemble de Cantor. On a en tête l'idée de montrer qu'il est de mesure nulle. On rappelle tout d'abord le lemme de régularité de Lebesgue qui va nous servir à formuler notre dictionnaire.

\vspace{3mm}

\paragraph*{\bf Premier lemme: points de densité} Si $A$ est un ensemble mesurable de l'espace euclidien $\mathbb{R}^n$, on dit que $x \in \mathbb{R}^n$ est un point de densité de $A$ si 

$$ \lim_{r \rightarrow 0}{\frac{\mathrm{Leb}\big( A \cap \mathrm{B}(x,r) \big)}{\mathrm{Leb}\big( \mathrm{B}(x,r) \big)}} \longrightarrow 1. $$  Le lemme suivant nous donne une information sur la structure locale d'un ensemble mesurable de $\mathbb{R}^n$ en terme de ses points de densité.

\begin{theorem}[Régularité de la mesure de Lebesgue]
Soit $A$ un ensemble de mesurable de $\mathbb{R}^n$. Alors presque tout point de $A$ est un point de densité de $A$. 
\end{theorem}

On donne maintenant un critère important qui lie les propriétés des orbites le long du flot de Teichmüller et l'ensemble $\mathcal{C}$.

\begin{lemma}
\label{densité}
Soit $\Sigma$ appartenant à $\mathcal{MD}_{g,n}$. Supposons que $\Theta(g_t \cdot \Sigma)$ ne tend pas vers $0$. Alors le point $\frac{\pi}{2}$ (représentant le feuilletage vertical de $\Sigma$) n'est pas un point de densité de $\mathcal{C}_{\Sigma} \subset S^1$. 

\end{lemma}

\begin{proof}

On donne ici une esquisse de preuve (mais qui ne cache aucune difficulté), une preuve rigoureuse avec $\epsilon$ peut-être trouvée dans \cite{Ghazman}[Lemma 19]. 

\vspace{2mm} \noindent On veut comprendre la structure locale de petits voisinages de $\frac{\pi}{2}$ dans $S^1$: le complémentaire de $\mathcal{C}_{\Sigma}$ est ouvert et ses composantes connexes correspondent à l'existence de cylindres dans les secteurs angulaires que ces composantes connexes représentent. 

\vspace{2mm} \noindent Le point $\frac{\pi}{2}$ est un point fixe répulsif de $g_t$ sur $\mathbb{RP}^1$ et donc son action nous permet de zoomer près de ce point; pour des $t$ très grands, les cylindres qu'on voit loin de la direction $\frac{\pi}{2}$ en étaient originellement très proches. L'hypothèse $\Theta(g_t \cdot \Sigma)$ ne tend pas vers $0$ nous assure que lorsqu'on zoome près de $\frac{\pi}{2}$, en "rembobinant" le flot, on voit des gros trous (par trous on entend composantes connexes du complémentaire de $\mathcal{C}_{\Sigma}$) près de $\frac{\pi}{2}$ ce qui l'empêche d'être un point de densité. 
\end{proof}

\vspace{2mm}

\paragraph*{\bf Deuxième lemme: directions hyperboliques}

On donne dans cette section un deuxième lemme, plus élémentaire que le premier, mais qui est crucial pour formuler l'analogie avec les variétés hyperboliques de volume infini qui est expliquée dans la section suivante.

\begin{lemma}
\label{vasque}
Soit $\Sigma \in \mathcal{MD}_{g,n}$ telle qu'il existe un cylindre affine en direction $\frac{\pi}{2}$ (\textit{i.e.} le feuilletage de direction $\frac{\pi}{2}$ contient une feuille fermée régulière hyperbolique). Alors $\Theta(g_t \cdot \Sigma)$ tend vers $\pi$ exponentiellement vite. 

\end{lemma} 

\begin{proof}
La preuve est laissée en exercice. Il s'agit juste de regarder ce que le flot de Teichmüller fait à un cylindre contenant la direction $\frac{\pi}{2}$; c'est un exercice d'algèbre linéaire élémentaire.
\end{proof}

\subsection{Analogie avec la géométrie hyperbolique}

Nous pensons qu'il y a une analogie intéressante entre les {\it espaces de modules de surfaces de dilatation munis du flot de Teichmüller} et les { \it variétés hyperboliques de volume infini munies de leur flot géodésique}\footnote{Pour être parfaitement rigoureux, il faudrait dire "les fibrés tangents de variétés des variétés de volume infini" car c'est bien sur cet espace que le flot géodésique agit. On s'affranchira de cette difficulté en disant juste "variétés hyperboliques"}.

\noindent Cette analogie peut-être justifiée et/ou motivée par les commentaires suivants:

\begin{enumerate}

\item il existe déjà une telle analogie pour les espaces de modules de surfaces de translation, mais avec des variétés de volume \textbf{fini};

\item le flot de Teichmüller restreint aux équivalents des surfaces de Veech dans le monde de dilatation (qu'on présente dans une section à venir) est \textit{exactement} un flot géodésique sur une surface hyperbolique géométriquement finie de volume infini;

\item le lemme \ref{vasque} met en place la première brique concrète de cette analogie: les surfaces dont le feuilletage vertical est Morse-Smale s'enfuit dans une "vasque" qui est formée des surfaces dont l'angle du plus grand cylindre est très grand.
\end{enumerate}

\noindent Une des motivations pour cette analogie est le théorème suivant (voir \cite{Ahlfors})

\begin{theorem}[Alfhors, \cite{Ahlfors}]
Soit $V = \mathbb{H}^3/\Gamma$ une variété hyperbolique de volume infini et de type géométriquement fini. Alors l'ensemble 
$$ \big\{ x \in T^1V \ | \ g_t \cdot x \ \text{ne pars pas à l'infini dans une vasque}  \big\}$$ est de mesure nulle.
\end{theorem}

\noindent En effet, si cette analogie était suffisamment robuste pour qu'un théorème analogue soit vrai dans le monde des surfaces de dilatation, nous pourrions surement prouver la généricité des feuilletages Morse-Smale (après laquelle, nous le rappelons au lecteur, nous courrons). En effet, comme l'indique le Lemme \ref{vasque}, on peut grossièrement identifier l'ensemble des feuilletage Morse-Smale à ceux dont l'orbite le long du flot de Teichmüller part à l'infini dans la "vasque" correspondant à l'angle. Il s'agit d'une heuristique qui doit être substantiellement raffinée pour être rigoureuse, mais qui est essentiellement juste. Le vrai problème est de montrer que le flot de Teichmüller agissant sur l'espace des modules de surfaces de dilatation ressemble structurellement suffisamment à un flot géodésique hyperbolique pour que le théorème ci-dessus se transporte d'un monde à l'autre. Nous illustrons dans la section suivante le seul exemple, en petite dimension, pour lequel nous savons faire fonctionner cette analogie.

\subsection{Le cas $g=1, n=2$}

On se restreint dans cette section au cas de tores avec deux points singuliers. On y discute la preuve du théorème suivant

\begin{theorem}
\label{morsesmalegenre1}
Pour presque tout tore de dilatation avec deux singularités, le flot vertical est Morse-Smale.
\end{theorem}

\noindent On ne donnera pas un preuve complète car le traitement de certains aspects techniques alourdirait considérablement l'exposition. Nous allons en discuter les grandes lignes. On réfère à \cite{Ghazman} pour une preuve précise et complète.

\vspace{3mm}

\paragraph*{\bf Une mesure invariante} Le point-clé qui nous permet de parvenir à une preuve de ce théorème est l'existence d'une mesure invariante par l'action de du flot de Teichmüller dans la classe de la mesure de Lebesgue qui rend l'analogie esquissée ci-dessus suffisamment robuste pour permettre un transport du théorème d'Alfhors. Quand $g=1$ et $n=2$, le feuilletage défini par l'action de $\mathrm{Sl}_2(\mathbb{R})$ sur $\mathcal{MD}_{1,2}(\lambda)$ coïncide avec le feuilletage par parties linéaires défini en \ref{partielineaire}. Parce que le feuilletage par parties linéaires est défini localement par des fonctions à valeurs dans le premier groupe de cohomologie de la surface associée, ce feuilletage porte une structure symplectique transverse. Couplé avec la mesure de Haar de $\mathrm{Sl}_2(\mathbb{R})$ poussée sur les feuilles, ce structure symplectique donne:

\begin{proposition}
Pour tout $\lambda > 1$ il existe sur $\mathcal{MD}_{1,2}(\lambda)$ une mesure $\mu$, invariante par l'action de $\mathrm{Sl}_2(\mathbb{R})$ et dans la même classe que la mesure de Lebesgue. De plus, la masse totale de $\mu$ est \textbf{infinie}.
\end{proposition}

\noindent Le fait que la masse de $\mu$ est infinie peut-être aisément vu de la manière suivante: l'action de $g_t$ est localement libre; en prenant un tore dont la direction verticale est Morse-Smale, il existe un petit ouvert $U$ autour de ce tore pour lequel tout élément a aussi une direction verticale Morse-Smale et donc est contenue dans un cylindre d'un certain angle variant continument dans $U$. Le flot de Teichmüller appliqué aux éléments de $U$ "gonfle" ce cylindre et repousse $U$ de plus en plus loin à l'infini (ceci est vu en utilisant la proposition \ref{vasque}, car la fonction $\Theta$ le long de $g_t$ va tendre vers $\pi$ uniformément pour les éléments de $U$). On peut ainsi s'assurer qu'il existe une suite de temps tendant vers l'infini $t_1, ..., t_n, ...$ tel que les $g_{t_i}(U)$ sont deux à deux disjoints. Comme $U$ est ouvert et que $\mu(U) > 0$, par $\sigma$-additivité la masse totale de $\mu$ doit être infinie.

\vspace{2mm}

\noindent Nous rappelons à notre lecteur que notre but est de démonter que $\mathcal{C}$ l'ensemble des points de $\mathcal{MD}_{1,2}(\lambda)$ dont le feuilletage vertical est minimal est de mesure zéro. 
La stratégie consiste à utiliser le lemme \ref{densité} couplé à l'existence de la mesure invariante $\mu$ pour montrer que l'ensemble des points de densité est de mesure nulle. Plus précisément on va se borner à montrer l'énoncé suivant

$$ \text{L'ensemble des points} \  \Sigma \in \mathcal{C} \ \text{qui sont des points de densité dans} \ \mathcal{C}_{\Sigma} \ \text{est de} \ \mu-\text{mesure nulle} $$ 

\noindent De cet énoncé le théorème se déduit aisément en utilisant Fubini. Le lemme \ref{densité} nous assure qu'un élément $\Sigma$ n'est pas un point de densité dans $\mathcal{C}_{\Sigma}$ dès lors que la trajectoire de $\Sigma$ via le flot de Teichmüller est suffisamment récurrente (précisément si la fonction $\Theta$ évaluée sur cette orbite ne tend pas vers $0$). Le point-clé est la proposition suivante 

\begin{proposition}

L'ensemble $$\mathcal{P} := \big\{ \Sigma \in \mathcal{MD}_{1,2}(\lambda) \ | \ \Theta(\Sigma) \leq \frac{\pi}{4} \big\}$$ est de $\mu$-mesure finie. 

\end{proposition}

 \noindent Nous ne donnerons pas de preuve de cette proposition qui est un calcul assez technique que le lecteur retrouvera dans \cite{Ghazman}. Nous expliquons son rôle dans l'analogie entre les  variétés hyperboliques et les espaces de modules de tores de dilatations. Elle assure que l'ensemble non-compact $\mathcal{P}$, dans lequel certaines trajectoires (pour le flot de Teichmüller) de surfaces dont le flot vertical est minimal peuvent s'enfuir, se comporte comme une pointe en géométrie hyperbolique, \textit{i.e.} est de volume fini. 

\vspace{3mm} \noindent Nous sommes maintenant en très bonne position pour prouver 

$$ \text{L'ensemble des points} \  \Sigma \in \mathcal{C} \ \text{qui sont des points de densité dans} \ \mathcal{C}_{\Sigma} \ \text{est de} \ \mu-\text{mesure nulle}.$$  En effet l'ensemble des $\Sigma$ qui sont des points de densité pour $\mathcal{C}_{\Sigma}$ doivent satisfaire $$\Theta(g_t \cdot \Sigma) \rightarrow 0.$$ Admettons que cet ensemble soit de mesure strictement positive. Alors le flot de Teichmüller devrait envoyer un ensemble de mesure strictement positive à l'infini dans $\mathcal{P}$. Comme ce dernier est de volume fini, ses intersections avec 

$$\big\{\Sigma \ | \ \Theta(\Sigma) \leq \epsilon\big\}$$ ont une mesure qui tend vers $0$ avec $\epsilon$. Il sera alors impossible de caser les images des points de densité par $g_t$ dans cet ensemble pour $t$ très grand, car $g_t$ préserve la mesure. Ceci contredit la possibilité que cet ensemble soit de mesure strictement positive. On a conclut la preuve du théorème \ref{morsesmalegenre1}.

\section{Action de $\mathrm{Sl}_2(\mathbb{R})$ et surfaces de Veech}

\label{secAct}
On discute dans cette section l'action de $\mathrm{Sl}_2(\mathbb{R})$ et ce que nous appelons les surfaces de Veech, qui sont une généralisation possible des surfaces de Veech au sens classique.

\subsection{Surfaces de Veech}

\begin{defi} 
Soit $\Sigma$ une surface de dilatation triangulable. On dit que $\Sigma$ est une surface de Veech si sa $\mathrm{Sl}_2(\mathbb{R})$-orbite est proprement plongée dans $\mathcal{MD}_{g,n}$. 
\end{defi}

\noindent Les surfaces de Veech sont donc les fermées invariants par $\mathrm{Sl}_2(\mathbb{R})$ de dimension minimale. Elles sont intéressantes pour deux raisons:

\begin{enumerate}

\item leur géométrie fournit un modèle simplifié pour la géométrie de l'espace de module entier;

\item le flot de Teichmüller s'y restreint en un flot géodésique sur une surface hyperbolique et est donc plus facile à comprendre.

\end{enumerate}

\noindent Nous ne connaissons à ce jour qu'assez peu de surfaces de Veech autres que celles de translation. Par ailleurs, toutes celles que nous connaissons sont des déformations de surfaces de Veech de translations. Nous posons donc les questions suivantes

\begin{question}
Existe-t-il des surfaces de Veech qui ne peuvent être déformées (via un chemin de surfaces de Veech) en une surface de Veech de translation?
\end{question}

\noindent De manière plus générale

\begin{question}

Existe-t-il des surfaces de Veech dont le groupe de Veech est trivial?

\end{question}

\noindent Ce que nous appelons le groupe de Veech d'une surface $\Sigma$ est son stabilisateur pour l'action de $\mathrm{Sl}_2(\mathbb{R})$. C'est un fait assez élémentaire que le groupe de Veech de toute surface triangulable est un sous-groupe discret de $\mathrm{Sl}_2(\mathbb{R})$, voir \cite{DFG}.

\vspace{3mm} \noindent  Nous donnons dans la suite de ce paragraphe trois exemples de surfaces de Veech. Nous ne prétendons pas que vérifier que ces surfaces sont effectivement des surfaces de Veech est évident, cela demande en générale une analyse de leur géométrie, dynamique et le calcul de certains éléments du groupe de Veech. Nous donnons cependant des références pour compléter les preuves.

\vspace{3mm}

\paragraph*{\bf Tores de Veech}  Nous donnons dans ce paragraphe un exemple de surfaces de Veech en genre $1$. Contrairement au cas de translation, grâce à la possibilité d'introduire dans singularité de type dilatation, il existe des surfaces et des espaces de modules intéressants dès le genre $1$. Les exemples qui suivent, dûs à Duryev, sont des tores de dilatation avec deux singularités dont le groupe de Veech est un sous groupe Zariski dense de $\mathrm{Sl}_2(\mathbb{R})$. Pour un paramètre $a\leq 1$, la surface $T_a$ est le tore de dilatation obtenu après recollement du l'hexagone ci-dessous

\begin{figure}[!h]
  \centering
  \includegraphics[scale=0.8]{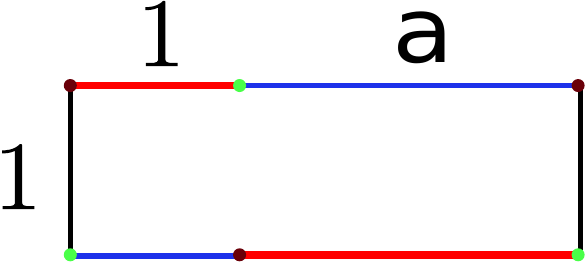}
  \caption{Un tore de Veech à deux singularités}
   \label{toreveech}
\end{figure}

Nous invitons le lecteur intéressé à vérifier les propriétés suivantes 

\begin{enumerate}

\item le groupe de Veech de $T_a$ est $$ \Gamma = \big\langle \begin{pmatrix}
1 & 1+a \\
0 & 1
\end{pmatrix}, \begin{pmatrix}
1 &  0\\
1 + \frac{1}{a} & 1
\end{pmatrix}, \begin{pmatrix}
-1 &  0\\
0 & -1
\end{pmatrix}  \big\rangle $$ 

\item La surface $T_a$ est une surface de Veech dont l'orbite via l'action de $\mathrm{Sl}_2(\mathbb{R})$(qui est fermée) est isomorphe à $\Gamma / \mathrm{Sl}_2(\mathbb{R})$.

\end{enumerate}

\noindent Le groupe $\Gamma$ est un groupe de co-volume infini de $\mathrm{Sl}_2(\mathbb{R})$. On conjecture que c'est toujours le cas pour une surface de Veech quelconque.

\begin{question}
Le groupe de Veech d'une surface de dilatation est-il toujours de covolume infini pour les surfaces qui ne sont pas de translation?
\end{question}

\noindent Il est remarquable que dans ce cas, le groupe de Veech permet d'expliciter intégralement la dynamique des feuilletages directionnels de $T_a$. Précisément, l'action de $\Gamma$ sur $S^1$ possède un \textit{ensemble limite} $\Lambda$ qui est son unique ensemble fermé minimal invariant. Cet ensemble  $\Lambda$ est un Cantor plongé dans $S^1$ et on a la trichotomie suivante

\begin{enumerate}

\item pour toute direction dans $S^1 \setminus \Lambda$, le feuilletage directionnel est Morse-Smale; 

\item pour toute direction dans le bord de $\Lambda$ (pas topologique, mais en tant que Cantor plongé dans $S^1$), le feuilletage à des connexions de selles fermées attractives et répulsives pour lesquelles le feuilletage est encore Morse-Smale;

\item pour les autres directions de $\Lambda$, le feuilletage est minimal.

\end{enumerate}

\noindent Cette caractérisation est une simple conséquence du fait que l'action du groupe de Veech sur l'ensemble des directions préserve les propriétés dynamiques. 

\vspace{3mm}

\paragraph*{\bf La surface à deux chambres identiques} On présente maintenant un type de surface de Veech qui est essentiellement différent de ce qu'on peut trouver dans le monde des surfaces de translation. C'est une surface de genre $2$ avec une singularité qui est obtenue par le collage suivant donné dans l'introduction (voir figure \ref{deuxchambres}).

\noindent Son groupe de Veech est engendré par l'élément parabolique suivant $\begin{pmatrix}
1&1\\
0&1
\end{pmatrix}$ et $\begin{pmatrix}
-1&0\\
0& -1
\end{pmatrix}$ la rotation d'angle $\pi$. Cet exemple contraste avec le cas des surfaces de translation: il s'agit d'une surface de Veech dont le groupe de Veech n'est pas Zariski dense. Il existe un Cantor de direction pour lesquel le feuilletage à un quasi-minimal qui est transversalement un Cantor. L'action du groupe de Veech est donc loin de "voir" toute la dynamique. On remarquera aussi qu'aucun des feuilletages directionnels n'est minimal. On réfère à \cite{DFG} et \cite{Lesguignols} pour plus de détails sur cet exemple.

\vspace{3mm}

\paragraph*{\bf La surface discothèque} 

La surface discothèque est probablement l'exemple le plus intéressant de surface de Veech que nous connaissons. Il mélange les deux types de caractéristiques observés avec les deux exemples précédents:

\begin{enumerate}

\item elles ont un gros groupe de Veech qui voit beaucoup de la dynamique: les directions correspondants à son ensemble limite correspondent à des feuilletages minimaux infiniment renormalisables ;

\item ce qui n'est pas vu par le groupe de Veech n'est pas complètement trivial: il peut s'agir soit de directions Morse-Smale soit de quasi-minimaux de genre intermédiaire. 

\end{enumerate}

\noindent La surface discothèque est la obtenue par le collage suivant

\begin{figure}[!h]
  \centering
  
  \includegraphics[scale=0.5]{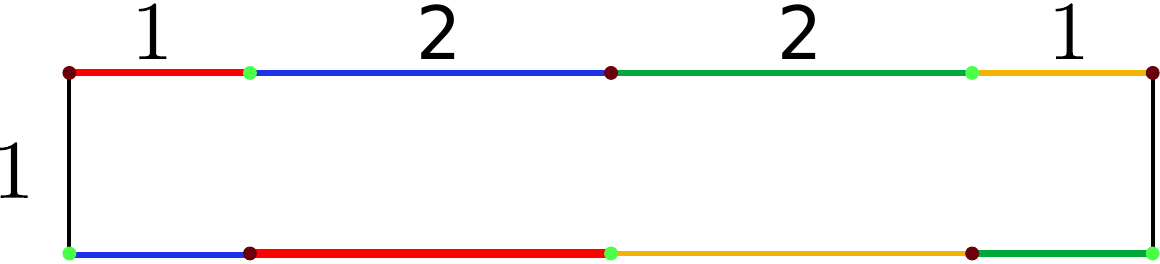}
  \caption{La fameuse discosurface}
   \label{discosurface}
\end{figure}

\noindent Son groupe de Veech est le groupe suivant 

$$ \big\langle \begin{pmatrix}
1 & 6 \\
0 & 1
\end{pmatrix}, \begin{pmatrix}
1 &  0\\
\frac{3}{2} & 1
\end{pmatrix}, \begin{pmatrix}
-1 &  0\\
0 & -1
\end{pmatrix}  \big\rangle. $$ C'est une surface de genre $2$ avec deux singularités de type euclidien.

\vspace{3mm}

\noindent Le groupe de Veech de la surface discothèque est un groupe de covolume infini dont l'ensemble limite est un ensemble de Cantor. Les directions de ce Cantor correspondent essentiellement à des directions minimales infiniment renormalisables. Les directions du complémentaires de cet ensemble limite sont presque toujours Morse-Smale, mais il existe un ensemble de dimension de Hausdorff nulle pour lesquelles les feuilletages directionnels s'accumulent sur des Cantor localisés dans des sous-surfaces de genre $1$. Nous renvoyons le lecteur à \cite{BFG} pour un traitement complet de cet exemple.

\subsection{Fermés $\mathrm{Sl}_2(\mathbb{R})$-invariants} De manière générale, l'action de $\mathrm{Sl}_2(\mathbb{R})$ nous semble être d'un intérêt certain. Nous connaissons finalement assez peu de choses sur les fermés invariants. Les deux seules familles d'exemples que nous connaissons sont assez élémentaires:

\begin{enumerate}

\item les orbites de surfaces de Veech;

\item les ensembles formés de surfaces dont le groupe d'holonomie linéaire est cyclique discret.

\end{enumerate}

\noindent Nous posons deux questions sur cette action.

\begin{question}
Pour toute composante connexe de $\mathcal{MD}_{g,n}(\lambda)$, existe-t-il toujours une $\mathrm{Sl}_2(\mathbb{R})$-orbite dense?
\end{question}

\noindent Nous connaissons au moins une famille d'exemples pour lesquels c'est le cas. Il s'agit de $\mathcal{MD}_{1,2}(\lambda)$ pour tout $\lambda$ pour lesquels toute feuille dont les surfaces n'ont pas une holonomie linéaire discrète est dense. Cette question devra peut-être reformulée pour éviter des contre-exemples triviaux en dimension plus grande (nous soupçonnons l'existence d'ouverts invariants stricts qui sont des raffinements des ouverts que forment les surfaces triangulables).

\begin{question}
Existe-t-il une $\mathrm{Sl}_2(\mathbb{R})$-orbite dont l'adhérence n'est pas une sous-variété lisse?
\end{question} 

\noindent Cette question est très intéressante, mais probablement un peu ambitieuse. On commence tout juste à étudier les fermés invariants pour les espaces homogènes de volume infini (voir \cite{McMullenMohammadiOh} par exemple).

\section{Lien avec la dynamique unidimensionnelle}
\label{secLien}

Nous tentons de présenter dans cette section, qui est essentiellement indépendante du reste du texte, les problématiques de dynamique unidimensionnelle auxquelles les surfaces de dilatations sont intimement liées. Nous supposons (soyons honnête) une connaissance basique du nombre de rotation et des échanges d'intervalles linéaires.

\vspace{3mm}

\noindent Les objets centraux de cette discussion sont les \textit{échanges d'intervalles généralisés} qui sont les bijections continues par morceaux de l'intervalle $[0,1]$ qui sont lisses sur leurs intervalles de continuité et qui préservent l'orientation. Ces systèmes dynamiques apparaissent naturellement comme applications de premier retour de flots sur des surfaces et sont donc des objets d'un grand intérêt. 

\begin{itemize}

\item La \textbf{théorie ergodique} de ses objets se ramène, via certaines considérations combinatoires assez élémentaires, à celle des \textit{échanges d'intervalles linéaires}. Ces objets ont été intensément étudiés au cours des dernières décennies, et sont au coeur du développement de la \textit{dynamique de Teichmüller}. 

\item En revanche, la \textbf{dynamique "lisse"} des échanges d'intervalles généralisés est moins bien comprise, et ne se réduit certainement pas aux échanges linéaires. Bien sûr, les propriétés ergodiques de ces objets interagissent fortement avec les questions de dynamique lisse mais ces dernières ne s'y réduisent pas. Les \textit{échanges d'intervalles affines} fournissent une famille assez simple pour appréhender cette problématique. 

\end{itemize}

\noindent En dynamique lisse, nous isolons deux questions qui nous semblent d'un intérêt tout particulier.

\begin{enumerate}

\item Comprendre le comportement dynamique \textit{générique} d'un tel objet.

\item Comprendre leurs classes de conjugaison $\mathcal{C}^1$.

\end{enumerate}

\subsection{Difféomorphismes du cercle}

Nous rappelons brièvement dans cette section la théorie des difféomorphismes du cercle, qui peuvent être pensés comme des échanges d'intervalles généralisés à deux branches. Tout d'abord, arrêtons nous sur les questions de théorie ergodique. La définition du \textbf{nombre de rotation} par Poincaré permet de régler la question ergodique: tout difféomorphisme de $S^1$ dont le nombre de rotation est irrationnel est semi-conjugué à la rotation de même nombre de rotation. Il est donc minimal et a une unique mesure invariante. 

\vspace{3mm} \noindent Intéressons nous donc à la question de la conjugaison $\mathcal{C}^1$. Partons d'un difféomorphisme $f : S^1 \longrightarrow S^1$ minimal et donc conjugué à $r_{\alpha} := x \mapsto x + \alpha \ [1]$ pour un certain $\alpha$ irrationnel. De prime abord, il peut sembler étrange que deux applications $\mathcal{C}^{\infty}$ puissent être conjuguée via autre chose qu'un difféomorphisme lisse. Cependant, cela à des implications assez fortes. En effet supposons 

$$f = \varphi \circ r_{\alpha} \circ \varphi^{-1}$$ avec $\varphi$ de classe au moins $\mathcal{C}^1$, il s'ensuit directement que les dérivées de $f^n$ sont uniformément bornée par $|| \varphi' || \cdot || (\varphi^{-1})' ||$ car $f^n = \varphi \circ r_{n\alpha} \circ \varphi^{-1}$. En appliquant simplement la règle de la chaine, on voit que $\log (f^n)'$ s'exprime comme les sommes de Birkhoff de $\log f'$ et la $\mathcal{C}^1$-conjugaison implique donc que les sommes de Birkhoff du logarithme de la dérivée sont bornée, ce qui est une restriction assez forte. 

\vspace{2mm}
\noindent Nous espérons avoir convaincu notre lecteur de la chose suivante: les premières obstructions à la conjugaison $\mathcal{C}^1$ sont à aller chercher dans la croissance des sommes de Birkhoff. Si $\mu$ est l'unique mesure invariante de $f$, il n'est pas difficile de se convaincre que systématiquement, dès lors que $f$ est minimal

$$ \int{\log f' d\mu} = 0. $$ On veut donc essayer de comprendre la croissance des sommes de Birkhoff de fonctions continues de moyenne nulle au-dessus d'une rotation irrationnelle.
\noindent Le théorème ergodique prédit que pour toute fonction $g$ continue de moyenne nulle et pour tout point $x$

$$ \sum_{i=0}^{n-1}{g(r_{\alpha}^i(x)} = o(n) $$ ce qui est une simple conséquence de l'unique ergodicité de $\mu$. Dans le cas des rotations, parce que le cercle est une variété très particulière, ce théorème peut être renforcé de la manière suivante 

\begin{lemma}[dit "de Denjoy-Koksma"]

Pour tout $\alpha$ il existe une suite de temps $0 < q_1 <\ldots < q_n < \ldots$ telle que pour toute fonction $g$ à variation bornée, de moyenne nulle on ait pour tout $n \in \mathbb{N}^*$

$$ \sum_{i=0}^{q_n}{g(r_{\alpha}^i(x)} \leq \mathrm{Var}(g).$$

\end{lemma} 

\noindent Autrement dit, il existe des temps spéciaux\footnote{Ces temps sont ceux pour lesquels $r_{\alpha}^{q_n}$ est très proche de l'identité} pour lesquels les sommes de Birkhoff de toute fonction suffisamment régulière sont bornées. C'est une propriété bien particulière qui nous dit qu'il n'y a pas de \textit{déviations} pour les sommes de Birkhoff. S'appuyant sur cette remarque, Herman développe dans le mémoire \textit{Sur la conjugaison différentiable des difféomorphismes du cercle à des rotations} une belle théorie qui à pour point culminant la preuve du théorème de conjugaison régulière suivant:

\begin{theorem}[Herman]

Pour $\alpha$ dans un ensemble de mesure pleine de $S^1$, tout difféomorphisme de classe $\mathcal{C}^{\infty}$ de nombre de rotation $\alpha$  est $\mathcal{C}^{\infty}$-conjugué à $r_{\alpha}$.

\end{theorem}

\noindent C'est un formidable théorème sur lequel Herman s'appuie pour prouver la chose suivante:

\begin{theorem}[Herman]

Soit $(f_t)_{t\in [0,1]}$ une famille à un paramètre de difféomorphisme du cercle variant de manière lisse telle que les nombres de rotation de $f_0$ et $f_1$ soient différents. Alors 

$$ \{ t \in [0,1] \ | \ f_t \text{est minimal} \} $$ est de mesure strictement positive. 

\end{theorem} 

\noindent Ce théorème prédit la coexistence de deux comportements dynamiques en probabilité: la minimalité et l'existence d'orbites périodiques.

\subsection{\'Echanges d'intervalles affines}

 Comme nous le disions plus haut, la théorie ergodique des échanges d'intervalles généralisés se réduit à celle des échanges linéaires, de manière analogue au fait que la théorie ergodique des difféomorphismes du cercle se réduit à celle des rotations. Le lecteur intéressé pourra consulter par exemple \cite{Yoccoz}. 

\vspace{3mm} \noindent Commençons par discuter la question de la régularité de la conjugaison. Soit $T$ un échange d'intervalle conjugué à un échange d'intervalle linéaire minimal $T_0$. Par le même raisonnement que dans le cas des rotations, les dérivées des itérés de $T$, les $T^n$, sont uniformément bornées. Cela nous invite une fois de plus à considérer la croissance des sommes de Birkhoff. Pour des raisons similaires au cas des difféomorphismes du cercle, le logarithme de la dérivée de $T$ est toujours d'intégrale nulle pour la mesure invariante. 

\noindent Dans ce cas, un programme initié par Zorich(voir \cite{Zorich2}) et complété par Forni(voir \cite{Forni}) a montré que les sommes de Birkhoff d'une fonction générique $h$ de moyenne nulle (pour un échange d'intervalle générique de genre $g \geq 2$) sont loin d'être bornées: il existe un nombre $0 <\alpha < 1$ tel que les sommes de Birkhoff de $h$ au temps $n$ sont de l'ordre de $n^{\alpha}$. Ceci suggère qu'il existe des obstructions de nature ergodique à la conjugaison régulière. 

\noindent La série d'articles \cite{CamelierGutierrez}, \cite{Cobo} \cite{BressaudHubertMaass}, \cite{MarmiMoussaYoccoz} montre que cette obstruction (qui vit dans la cohomologie de la surface associée) est réalisée dans la classe des échanges d'intervalles affines. C'est une indication que cette classe voit déjà une partie de la richesse de la théorie générale. Les surfaces de dilatation, pour lesquelles les applications de premier retour des flots directionnels sont des échanges d'intervalles affines, donnent une réalisation géométrique des ces systèmes dynamiques. La question analogue à celle qu'on pose pour les feuilletages de surfaces de dilatation est la suivante:

\begin{question}

Dans la famille des échanges d'intervalles affine, le sous-ensemble formés des échanges de dynamique Morse-Smale est-il de mesure pleine?

\end{question}

\noindent Pour donner une dernière raison qui nous semble justifier la place centrale des échanges d'intervalles affines dans la théorie plus générale des échanges généralisés, notons qu'ils sont de bons candidats pour être un attracteur de l'opérateur de renormalisation, à savoir de l'induction de Rauzy. Ce fait a été prouvé dans le cas des combinatoires de genre $1$, voir \cite{CunhaSmania}.

\bibliographystyle{alpha} 
\bibliography{biblio}

\end{document}